\documentclass{gtart}

\def\ifplaintex{\expandafter\ifx\csname documentclass\endcsname\relax}

\ifplaintex 
\else
\fi

\expandafter\ifx\csname beginpicture\endcsname\relax
\expandafter\ifx\csname documentclass\endcsname\relax
\input pictex \else
\input prepictex \input pictex \input postpictex \fi\fi

\def\gt{{\mathsurround=0pt\it $\cal G\mskip-2mu$eometry \&\ 
$\cal T\!\!$opology}}        

\def\gtp{{\mathsurround=0pt\it $\cal G\mskip-2mu$eometry \&\ 
$\cal T\!\!$opology $\cal P\!$ublications}}  

\def\lognumber#1{\def\thelognumber{#1}}
\def\volumenumber#1{\def\thevolumenumber{#1}}
\def\papernumber#1{\def\thepapernumber{#1}}
\def\volumeyear#1{\def\thevolumeyear{#1}}

\def\pagenumbers#1#2{\def\startpage{#1}\def\finishpage{#2}}
\def\published#1{\def\publishdate{#1}}
\def\proposed#1{\def\theproposer{#1}}
\def\seconded#1{\def\theseconders{#1}}
\def\received#1{\def\receiveddate{#1}}
\def\revised#1{\def\reviseddate{#1}}
\def\accepted#1{\def\accepteddate{#1}}

\def\asciiaddress#1{\def\theasciiaddress{#1}}

\long\def\asciiabstract#1{\long\def\theasciiabstract{#1}}

\let\\\par\let\thelognumber\relax
\let\thevolumenumber\relax\let\thepapernumber\relax
\let\thevolumeyear\relax\let\thesamplenumber\relax\let\startpage\relax
\let\finishpage\relax\let\publishdate\relax\let\receiveddate\relax
\let\reviseddate\relax\let\accepteddate\relax\let\theasciititle\relax
\let\theasciiauthors\relax\let\theasciiaddress\relax
\let\theasciiabstract\relax
\let\theasciiemail\relax\let\theshortauthors\relax\let\theshorttitle\relax

\long\def\maketitlep{   

\count0=\startpage

\gt\hfill      
\beginpicture
\setcoordinatesystem units <0.33truein, 0.33truein> point at 2.2 0.9
\setplotsymbol ({$\cal G$})
\plotsymbolspacing=9truept
\circulararc 315 degrees from 0 1 center at 0 0
\setplotsymbol ({$\cal T$})
\circulararc 315 degrees from 1 -1 center at 1 0
\endpicture
\break
{\small\ifx\thesamplenumber\relax 
Volume \else Sample
\fi\thevolumenumber\ (\thevolumeyear)
\startpage--\finishpage\nl
Published: \publishdate}
\vglue 0.5truein plus 0.4fil minus 0.1truein

{\parskip=0pt\leftskip 0pt plus 1fil\def\\{\par\smallskip}{\ifplaintex\large
\else\Large\fi\bf\thetitle}\par\medskip}   

\vglue 0pt plus 0.1fil 

{\parskip=0pt\leftskip 0pt plus 1fil\def\\{\par}{\sc\theauthors}
\par\medskip}

\vglue 0pt plus 0.1fil 

{\small\parskip=0pt\let\newline\\
{\leftskip 0pt plus 1fil\def\\{\par}{\sl\theaddress}\par}
\expandafter\ifx\theemail\relax    
\relax\else\vglue 5pt plus 0.02fil minus 2pt\def\\{\stdspace{\rm 
and}\stdspace} 
\cl{Email:\stdspace\tt\theemail}\fi
\ifx\theurl\relax                  
\relax\else\vglue 5pt plus 0.02fil minus 2pt\def\\{\stdspace{\rm 
and}\stdspace}
\cl{URL:\stdspace\tt\theurl}\fi\par}

\vglue 7pt plus 0.3fil minus 3pt

{\bf Abstract}
\vglue 5pt plus 0.1fil minus 2pt

\theabstract

\vglue 7pt plus 0.3fil minus 3pt

{\bf AMS Classification numbers}\quad Primary:\quad \theprimaryclass

Secondary:\quad \thesecondaryclass

\vglue 5pt plus 0.3fil minus 2pt

{\bf Keywords}\quad \thekeywords

\vglue 10pt plus 0.5fil minus 5pt

{\small  Proposed: \theproposer\hfill Received: \receiveddate\nl
Seconded: \theseconders\hfill 
\ifx\reviseddate\relax                         
Accepted: \accepteddate                        
\else
Revised: \reviseddate                          
\fi}
\eject
}       

\let\maketitlepage\maketitlep
\let\maketitle\maketitlepage

\font\phead=cmsl9 scaled 950
\font=cmsl9 scaled 1050
\font\pnum=cmbx10 scaled 913
\font\lnum=cmbx10 
\font\pfoot=cmsl9 scaled 950
\font=cmsl9 scaled 1050
\ifplaintex
\headline{\vbox to 0pt{\vskip -4.5mm\line{\small\phead\ifnum
\count0=\startpage ISSN 1364-0380 (on line)
1465-3060 (printed) \hfill {\pnum\folio}\else\ifodd\count0\def\\{ }%
\ifx\theshorttitle\relax\thetitle\else\theshorttitle\fi\hfill{\pnum\folio}
\else\def\\{ and }{\pnum\folio}\hfill\ifx\theshortauthors\relax\theauthors
\else\theshortauthors\fi\fi\fi}\vss}}
\footline{\vbox to 0pt{\vglue 0mm\line{\small\pfoot\ifnum\count0=\startpage
\copyright\ \gtp\hfill\else
\gt, Volume \thevolumenumber\ (\thevolumeyear)\hfill\fi}\vss
}}
\else
\makeatletter
\def\@oddhead{{\small\lhead\ifnum\count0=\startpage ISSN 1364-0380 (on line)
1465-3060 (printed) \hfill {\lnum\number\count0}\else\ifodd\count0
\def\\{ }\ifx\theshorttitle\relax \thetitle \else\theshorttitle\fi\hfill
{\lnum\number\count0}\else\def\\{ and }{\lnum\number\count0}
\hfill\ifx\theshortauthors\relax 
\theauthors\else\theshortauthors\fi\fi\fi}}\def\@evenhead{\@oddhead}
\def\@oddfoot{\small\lfoot\ifnum\count0=\startpage\copyright\ \gtp\hfill\else
\gt, Volume \thevolumenumber\ (\thevolumeyear)\hfill\fi}
\def\@evenfoot{\@oddfoot}
\makeatother
\fi

\newwrite\gtoutfile
\long\gdef\makeheadfile{  
{\def\\{, }\def\s{ }
\immediate\openout\gtoutfile head.xxx
\immediate\write\gtoutfile{To: math@arxiv.org}
\immediate\write\gtoutfile{Subject: put or rep NNNNN:pppp}
\immediate\write\gtoutfile{--text follows this line--}
\immediate\write\gtoutfile{Proxy-for: \ifx\theasciiauthors\relax
\theauthors\else\theasciiauthors\fi\s<\ifx\theasciiemail\relax\theemail\else\theasciiemail\fi>}
\immediate\write\gtoutfile{\noexpand\\}
\immediate\write\gtoutfile{Authors: \ifx\theasciiauthors\relax
\theauthors\else\theasciiauthors\fi}
{\def\\{ }\immediate\write\gtoutfile{Title: \ifx\theasciititle\relax
\thetitle\else\theasciititle\fi}}
\immediate\write\gtoutfile{Subj-class: GT or SG or MG etc}
\immediate\write\gtoutfile{MSC-class: \theprimaryclass\ifx\thesecondaryclass\relax\else, \thesecondaryclass\fi}
\immediate\write\gtoutfile{Journal-ref: Geom. Topol. \thevolumenumber
(\thevolumeyear) \startpage-\finishpage}
\immediate\write\gtoutfile{Comments: Published by Geometry and Topology at}
\immediate\write\gtoutfile{\s\s http://www.maths.warwick.ac.uk/gt/GTVol\thevolumenumber/paper\thepapernumber.abs.html}
\immediate\write\gtoutfile{\noexpand\\}
\immediate\write\gtoutfile{}
\ifx\theasciiabstract\relax
\immediate\write\gtoutfile{\theabstract}\else
\immediate\write\gtoutfile{\theasciiabstract}\fi
\immediate\write\gtoutfile{}
\immediate\write\gtoutfile{\noexpand\\}
\immediate\write\gtoutfile{}
\immediate\closeout\gtoutfile}}  

\def\maketitlepage{\maketitlep\makeheadfile}
\let\maketitle\maketitlepage

\def\ifplaintex{\expandafter\ifx\csname documentclass\endcsname\relax}

\ifplaintex 
\else
\fi

\expandafter\ifx\csname beginpicture\endcsname\relax
\expandafter\ifx\csname documentclass\endcsname\relax
\input pictex \else
\input prepictex \input pictex \input postpictex \fi\fi

\def\gt{{\mathsurround=0pt\it $\cal G\mskip-2mu$eometry \&\ 
$\cal T\!\!$opology}}        

\def\gtp{{\mathsurround=0pt\it $\cal G\mskip-2mu$eometry \&\ 
$\cal T\!\!$opology $\cal P\!$ublications}}  

\def\lognumber#1{\def\thelognumber{#1}}
\def\volumenumber#1{\def\thevolumenumber{#1}}
\def\papernumber#1{\def\thepapernumber{#1}}
\def\volumeyear#1{\def\thevolumeyear{#1}}

\def\pagenumbers#1#2{\def\startpage{#1}\def\finishpage{#2}}
\def\published#1{\def\publishdate{#1}}
\def\proposed#1{\def\theproposer{#1}}
\def\seconded#1{\def\theseconders{#1}}
\def\received#1{\def\receiveddate{#1}}
\def\revised#1{\def\reviseddate{#1}}
\def\accepted#1{\def\accepteddate{#1}}

\def\asciiaddress#1{\def\theasciiaddress{#1}}

\long\def\asciiabstract#1{\long\def\theasciiabstract{#1}}

\let\\\par\let\thelognumber\relax
\let\thevolumenumber\relax\let\thepapernumber\relax
\let\thevolumeyear\relax\let\thesamplenumber\relax\let\startpage\relax
\let\finishpage\relax\let\publishdate\relax\let\receiveddate\relax
\let\reviseddate\relax\let\accepteddate\relax\let\theasciititle\relax
\let\theasciiauthors\relax\let\theasciiaddress\relax
\let\theasciiabstract\relax
\let\theasciiemail\relax\let\theshortauthors\relax\let\theshorttitle\relax

\long\def\maketitlep{   

\count0=\startpage

\gt\hfill      
\beginpicture
\setcoordinatesystem units <0.33truein, 0.33truein> point at 2.2 0.9
\setplotsymbol ({$\cal G$})
\plotsymbolspacing=9truept
\circulararc 315 degrees from 0 1 center at 0 0
\setplotsymbol ({$\cal T$})
\circulararc 315 degrees from 1 -1 center at 1 0
\endpicture
\break
{\small\ifx\thesamplenumber\relax 
Volume \else Sample
\fi\thevolumenumber\ (\thevolumeyear)
\startpage--\finishpage\nl
Published: \publishdate}
\vglue 0.5truein plus 0.4fil minus 0.1truein

{\parskip=0pt\leftskip 0pt plus 1fil\def\\{\par\smallskip}{\ifplaintex\large
\else\Large\fi\bf\thetitle}\par\medskip}   

\vglue 0pt plus 0.1fil 

{\parskip=0pt\leftskip 0pt plus 1fil\def\\{\par}{\sc\theauthors}
\par\medskip}

\vglue 0pt plus 0.1fil 

{\small\parskip=0pt\let\newline\\
{\leftskip 0pt plus 1fil\def\\{\par}{\sl\theaddress}\par}
\expandafter\ifx\theemail\relax    
\relax\else\vglue 5pt plus 0.02fil minus 2pt\def\\{\stdspace{\rm 
and}\stdspace} 
\cl{Email:\stdspace\tt\theemail}\fi
\ifx\theurl\relax                  
\relax\else\vglue 5pt plus 0.02fil minus 2pt\def\\{\stdspace{\rm 
and}\stdspace}
\cl{URL:\stdspace\tt\theurl}\fi\par}

\vglue 7pt plus 0.3fil minus 3pt

{\bf Abstract}
\vglue 5pt plus 0.1fil minus 2pt

\theabstract

\vglue 7pt plus 0.3fil minus 3pt

{\bf AMS Classification numbers}\quad Primary:\quad \theprimaryclass

Secondary:\quad \thesecondaryclass

\vglue 5pt plus 0.3fil minus 2pt

{\bf Keywords}\quad \thekeywords

\vglue 10pt plus 0.5fil minus 5pt

{\small  Proposed: \theproposer\hfill Received: \receiveddate\nl
Seconded: \theseconders\hfill 
\ifx\reviseddate\relax                         
Accepted: \accepteddate                        
\else
Revised: \reviseddate                          
\fi}
\eject
}       

\let\maketitlepage\maketitlep
\let\maketitle\maketitlepage

\font\phead=cmsl9 scaled 950
\font=cmsl9 scaled 1050
\font\pnum=cmbx10 scaled 913
\font\lnum=cmbx10 
\font\pfoot=cmsl9 scaled 950
\font=cmsl9 scaled 1050
\ifplaintex
\headline{\vbox to 0pt{\vskip -4.5mm\line{\small\phead\ifnum
\count0=\startpage ISSN 1364-0380 (on line)
1465-3060 (printed) \hfill {\pnum\folio}\else\ifodd\count0\def\\{ }%
\ifx\theshorttitle\relax\thetitle\else\theshorttitle\fi\hfill{\pnum\folio}
\else\def\\{ and }{\pnum\folio}\hfill\ifx\theshortauthors\relax\theauthors
\else\theshortauthors\fi\fi\fi}\vss}}
\footline{\vbox to 0pt{\vglue 0mm\line{\small\pfoot\ifnum\count0=\startpage
\copyright\ \gtp\hfill\else
\gt, Volume \thevolumenumber\ (\thevolumeyear)\hfill\fi}\vss
}}
\else
\makeatletter
\def\@oddhead{{\small\lhead\ifnum\count0=\startpage ISSN 1364-0380 (on line)
1465-3060 (printed) \hfill {\lnum\number\count0}\else\ifodd\count0
\def\\{ }\ifx\theshorttitle\relax \thetitle \else\theshorttitle\fi\hfill
{\lnum\number\count0}\else\def\\{ and }{\lnum\number\count0}
\hfill\ifx\theshortauthors\relax 
\theauthors\else\theshortauthors\fi\fi\fi}}\def\@evenhead{\@oddhead}
\def\@oddfoot{\small\lfoot\ifnum\count0=\startpage\copyright\ \gtp\hfill\else
\gt, Volume \thevolumenumber\ (\thevolumeyear)\hfill\fi}
\def\@evenfoot{\@oddfoot}
\makeatother
\fi

\newwrite\gtoutfile
\long\gdef\makeheadfile{  
{\def\\{, }\def\s{ }
\immediate\openout\gtoutfile head.xxx
\immediate\write\gtoutfile{To: math@arxiv.org}
\immediate\write\gtoutfile{Subject: put or rep NNNNN:pppp}
\immediate\write\gtoutfile{--text follows this line--}
\immediate\write\gtoutfile{Proxy-for: \ifx\theasciiauthors\relax
\theauthors\else\theasciiauthors\fi\s<\ifx\theasciiemail\relax\theemail\else\theasciiemail\fi>}
\immediate\write\gtoutfile{\noexpand\\}
\immediate\write\gtoutfile{Authors: \ifx\theasciiauthors\relax
\theauthors\else\theasciiauthors\fi}
{\def\\{ }\immediate\write\gtoutfile{Title: \ifx\theasciititle\relax
\thetitle\else\theasciititle\fi}}
\immediate\write\gtoutfile{Subj-class: GT or SG or MG etc}
\immediate\write\gtoutfile{MSC-class: \theprimaryclass\ifx\thesecondaryclass\relax\else, \thesecondaryclass\fi}
\immediate\write\gtoutfile{Journal-ref: Geom. Topol. \thevolumenumber
(\thevolumeyear) \startpage-\finishpage}
\immediate\write\gtoutfile{Comments: Published by Geometry and Topology at}
\immediate\write\gtoutfile{\s\s http://www.maths.warwick.ac.uk/gt/GTVol\thevolumenumber/paper\thepapernumber.abs.html}
\immediate\write\gtoutfile{\noexpand\\}
\immediate\write\gtoutfile{}
\ifx\theasciiabstract\relax
\immediate\write\gtoutfile{\theabstract}\else
\immediate\write\gtoutfile{\theasciiabstract}\fi
\immediate\write\gtoutfile{}
\immediate\write\gtoutfile{\noexpand\\}
\immediate\write\gtoutfile{}
\immediate\closeout\gtoutfile}}  

\def\maketitlepage{\maketitlep\makeheadfile}
\let\maketitle\maketitlepage

\lognumber{197}

\volumenumber{6}
\papernumber{24} 
\volumeyear{2002}
\pagenumbers{815}{852} 
\received{16 July 2001}
\revised{9 December 2002}
\accepted{18 December 2002}
\published{18 December 2002}
\proposed{Jean-Pierre Otal}
\seconded{David Gabai, Walter Neumann}

\usepackage{epsf}
\usepackage{amsmath} 
\usepackage{amssymb}
\usepackage{amscd}

\usepackage{enumerate}

\newtheorem{thm}{Theorem}[section]    
\newtheorem{lem}[thm]{Lemma}         
\newtheorem{prop}[thm]{Proposition}
\newtheorem{cor}[thm]{Corollary}
\newtheorem*{theorema}{Theorem A}   
\newtheorem*{theoremb}{Theorem B}

\newtheorem*{adtheoremb}{Addendum to Theorem B}
\theoremstyle{definition}
\newtheorem{defn}[thm]{Definition}    
\newtheorem*{rem}{Remark}             %
\newtheorem*{conv}{Convention}

\def\VERT{\operatorname{VERT}}
\def\TRANS{\operatorname{TRANS}}
\def\Isom{\operatorname{Isom}}
\def\ROT{\operatorname{ROT}}

\def\hol{\operatorname{hol}}

\def\NN{\mathop{\mathcal{N}}}
\def\Hom{\operatorname{Hom}}

\def\im{\operatorname{Im}}
\def\re{\operatorname{Re}}
\def\OO{\mathcal{O}}

\def\C{{\mathbb C}}

\def\0{\emptyset}
\def\R{{\mathbb R}}

\def\H{{\mathbb H}}
\def\Z{{\mathbb Z}}

\def\Q{{\mathbb Q}}

\def\3{\ss}
\def\8{\infty}

\def\<{\langle}
\def\>{\rangle}

\begin{document}
\title{Regenerating hyperbolic cone structures from Nil}
\author{Joan Porti}
\address{Departament de Matem\`atiques,  Universitat
Aut\`onoma de Barcelona\\08193 Bellaterra, Spain}
\asciiaddress{Departament de Matematiques,  Universitat
Autonoma de Barcelona\\08193 Bellaterra, Spain}

\email{porti@mat.uab.es}

\begin{abstract}
Let $\mathcal O$ be a three-dimensional $Nil$--orbifold, with
branching locus a knot $\Sigma$ transverse to the Seifert fibration.
We prove that $\mathcal O$ is the limit of hyperbolic cone manifolds
with cone angle in $(\pi-\varepsilon,\pi)$.  We also study the space
of Dehn filling parameters of $\mathcal O-\Sigma$. Surprisingly it is
not diffeomorphic to the deformation space constructed from the
variety of representations of $\mathcal O-\Sigma$.
 As a corollary of this,  we find examples of
spherical cone manifolds with singular set a knot that  are not
locally rigid. Those examples have large cone angles.
\end{abstract}

\asciiabstract{Let O be a three-dimensional Nil-orbifold, with
branching locus a knot Sigma transverse to the Seifert fibration.  We
prove that O is the limit of hyperbolic cone manifolds with cone angle
in (pi-epsilon, pi).  We also study the space of Dehn filling
parameters of O-Sigma.  Surprisingly it is not diffeomorphic to the
deformation space constructed from the variety of representations of
O-Sigma.  As a corollary of this, we find examples of spherical cone
manifolds with singular set a knot that are not locally rigid. Those
examples have large cone angles.}



\maketitlepage

\section{Introduction}

This paper is motivated by a phenomenon occurring in the proof of the
orbifold theorem. This proof suggests that some orbifolds with
geometry $Nil$ appear as limit of rescaled hyperbolic cone manifolds.
In the current proofs of the orbifold theorem
\cite{BLP0,BLP1,BLP2,BP,CHK}, it is only shown that those families of
cone manifolds collapse, and this is used to construct a Seifert
fibration of the orbifold, without knowing which kind of geometric
structure is involved.

Every closed three-dimensional $Nil$ orbifold admits an orbifold
Seifert fibration. We assume that the ramification locus is a circle
transverse to its Seifert fibration. This implies that the
ramification index is 2. Hence we view the orbifold as a cone manifold
with cone angle $\pi$.

\begin{theorema}\label{thm:a}
Let $\mathcal O$ be a closed three-dimensional $Nil$ orbifold
whose ramification locus $\Sigma$ is a circle transverse to its
Seifert fibration. Then there exist a family of hyperbolic cone
structures on the underlying space of $\mathcal O$
 with singular set $\Sigma$ parametrized by the cone angle
$\alpha\in(\pi-\varepsilon,\pi)$, for some $\varepsilon>0$.

In addition, when $\alpha\to \pi^-$ these hyperbolic cone
manifolds converge to a point. If they are re-scaled by
$(\pi-\alpha)^{-1/3}$, then they converge to a Euclidean
2--orbifold, which is the basis of the Seifert fibration of
$\mathcal O$. Finally, if they are re-scaled by
$(\pi-\alpha)^{-1/3}$ in the horizontal direction and
$(\pi-\alpha)^{-2/3}$ in the vertical one, then they converge to
$\OO$.
\end{theorema}

If the ramification locus $\Sigma $ was a circle but not transverse
to the Seifert fibration of $\OO$, then $\Sigma$ would be a fibre.
 In this case the conclusion of Theorem~A could not hold, because
  $\OO-\Sigma$ must be
hyperbolic, and therefore $\OO-\Sigma$ can not be Seifert fibred.

The following corollary follows from Theorem~A and Kojima's global
rigidity theorem \cite{Kojima}.

\begin{cor}
Let $\mathcal O$ be an orbifold as in Theorem~A.
 There exist a family of hyperbolic cone structures on the
underlying space of $\mathcal O$
 with singular set $\Sigma$ parametrized by the cone angle
$\alpha\in(0,\pi)$.
\end{cor}

The first part of Theorem~A is a particular case of  Theorem~B below,
which gives a larger space of deformations parametrized by
Dehn-filling coefficients. A cone manifold structure on
$\vert\mathcal O\vert$ with singular set $\Sigma$ induces a
non-complete metric on $\mathcal O-\Sigma$, whose completion is
precisely the cone manifold. This is a particular case of
structures on the end of $\mathcal O-\Sigma$ called of \emph{Dehn
type}. Those structures are defined by Thurston in
\cite{ThurstonNotes} and they are described by a pair
$(p,q)\in\mathbb R^2\cup\{\infty\}$.

\begin{theoremb} Let $\mathcal O$ be a $Nil$ 3--orbifold as in
Theorem~A. There exists a neighborhood $U$ of $(2,0)$  in $\mathbb
R^2$ and two $\mathcal C^1$--functions $f\co
(-\varepsilon,\varepsilon)\to (-\infty,2]$ concave and $g\co
(-\varepsilon,\varepsilon)\to [2,+\infty)$ convex, with
$f\vert_{[0,\varepsilon)} \equiv g\vert_{[0,\varepsilon)} \equiv 2$
and
\[
\lim \limits_{q\to
    0^-} \frac{2-f(q)} {\vert q \vert ^{3/2}}
    =
    \lim \limits_{q\to
    0^-} \frac{g(q)-2}{\vert q \vert ^{3/2}}>0,
\]
such that the following hold. Every point in $\{(p,q)\in U\mid p\geq
f(q)\}$ is the Dehn-filling coefficient of a
 geometric structure on $\mathcal O-\Sigma$ of the following kind:

\begin{itemize}
\item[-] hyperbolic for $p>g(q)$;
\item[-] Euclidean for $p=g(q)$, $q<0$;
\item[-] spherical for $p<g(q)$, $q>0$.
\end{itemize}
In addition, every point in the line $p=2$ corresponds to a
transversely Riemannian foliation of  codimension two
 (transversely hyperbolic for $q>0$, Euclidean for $q=0$ and spherical for
 $q>0$).
\end{theoremb}

 \begin{figure}[htbp]
 \epsfxsize 2.3in
 \centerline{\epsfbox{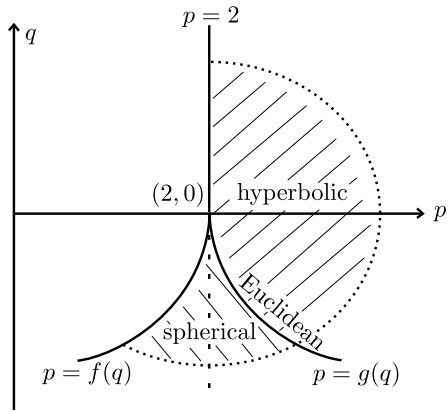}}
 \caption{The open set of Theorem~B} \protect \label{fig:dehnf}
 \end{figure}

When $q=0$, Dehn filling coefficients $(p,0)$ correspond to cone
structures with cone angle $2\pi/p$. Hence Theorem~B implies the
existence of hyperbolic cone manifolds with cone angles in
$(\pi-\varepsilon,\pi)$ of Theorem~A.

\begin{figure}[htbp]
 \epsfxsize.8\hsize
 \centerline{\epsfbox{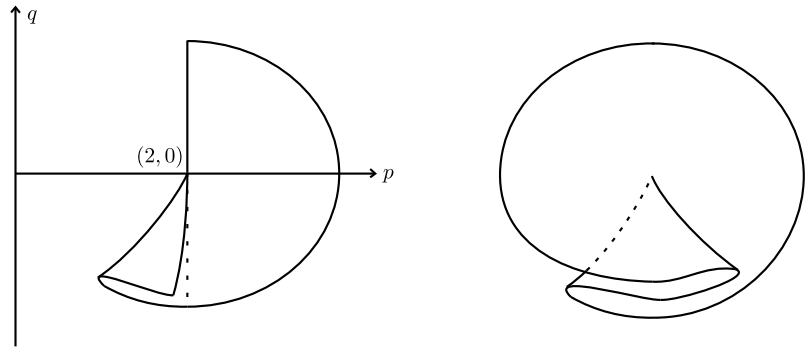}} \caption{
 The picture on the right hand side represents a Whitney pleat (a map conjugate to
 $(x,y)\mapsto (x,y^3-xy)$). The picture on the left
 hand side
 shows the situation in Theorem~B: half of it.} \protect \label{fig:Whitney}
 \end{figure}

To prove Theorem~B, we construct a deformation space
 homeomorphic to
a half-disc. However, Dehn filling coefficients do not define a
homeomorphism between the deformation space and the region of
Theorem~B, because there is a Whitney pleat at the point
$(p,q)=(2,0)$ corresponding to the $Nil$ structure (see
Figure~\ref{fig:Whitney}).

The image of the folding region is precisely the the curve $p=f(q)$,
$ q<0$. Thus we have the following addendum to Theorem~B:

\begin{adtheoremb} Local rigidity fails to hold on the curve $p=f(q)$,
$ q<0$. In addition, every Dehn filling coefficient in
\[\{(p,q)\in U\mid f(q)<p<2\}\] corresponds to two different spherical structures,
and every Dehn filling coefficient $(2,q)$ with $q<0$ corresponds to
a spherical structure and a transversely spherical foliation.
\end{adtheoremb}

By considering straight lines with rational slopes that intersect the
curve $ p=f(q)$, $q<0$ we obtain the following corollary.

\begin{cor} Local rigidity fails to hold for some spherical cone manifolds
with singular set a knot and large cone angles.
\end{cor}

In 1998 Casson  showed that local rigidity fails for some
hyperbolic cone manifolds with singular set a graph. Local
rigidity for compact hyperbolic cone manifolds with singular set a
link and cone angles $\leq 2\pi$ has been proved by Hodgson and
Kerckhoff in \cite{HK}. Likely, their methods can be adapted to
the situation in the spherical case, but our corollary shows that
an upper bound of the cone angle is essential in the spherical
case.

The proof of Theorems~A and~B allows to prove the the following
metric properties of the family of collapsing cone manifolds.

 \begin{prop}\label{prop:properties1}
 Let $C_{\alpha}$ denote the hyperbolic cone manifolds provided by
 Theorem~A, with $\alpha\in(\pi-\varepsilon,\pi)$, and let
 $\Sigma_{\alpha}$ denote its singular set. Then:
\[
\lim_{\alpha\to\pi^-} \frac{\operatorname{vol}(C_{\alpha})}{
(\pi-\alpha)\operatorname{length}(\Sigma_{\alpha})}=\frac38
\quad\textrm{ and }\quad \lim_{\alpha\to\pi^-}
\frac{\operatorname{length}(\Sigma_{\alpha})}{
(\pi-\alpha)^{1/3}}=l_0>0.
\]
 \end{prop}

\begin{prop}\label{prop:properties2}
 Let $f$ and $g$ be the functions of Theorem~B and let $l_0>0$ be as in previous
 proposition. Then:
\[
 \lim \limits_{q\to
    0^-}  \frac{2-f(q)}{\vert q \vert ^{3/2}}
=
 \lim \limits_{q\to
    0^-}  \frac{g(q)-2}{\vert q \vert ^{3/2}}
    =\frac {4}{9\sqrt[4]
3\pi}\,
    l_0^{3/2}.
\]
\end{prop}

In the proof of Theorem~B, we   first construct spaces of
geometric structures on $M$ parametred by $(s,t)\in U\subset
\mathbb R^2$, where $U$ a neighborhood of the origin
(Theorem~\ref{thm:defstr}). We consider spaces of both, hyperbolic
and spherical structures, and work with unified notation: $\mathbb
X^3$   denotes either $\mathbb H^3$ or $\mathbb S^3$. Those
structures are non degenerate except when $s=0$ or $t=0$. The
degenerated  structures  are the following ones: the origin
corresponds to the original $Nil$ structure, the line $t=0$ to
Euclidean structures, and the line $s=0$ to transversely
hyperbolic or spherical foliations. This space of structures $U$
has symmetry: $(\pm s,\pm t)$ is the parameter of the same
structure as $(s,t)$ up to changing the orientation or the spin
structure.

Next we construct a deformation space {\it Def}, which is a half disc
centered at the origin with parameters $(s,\tau)$, with $s\geq 0$,
and $\tau=t^2$ in the hyperbolic case, $\tau=-t^2$ in the
spherical case, and $\tau=0$ in the Euclidean one. In the proof of
Theorem~B, we show that the Dehn filling coefficients $(p,q)$
define an analytic map on $(s,\tau)$ that has ``half Whitney
pleat" at the origin, as illustrated in Figure~\ref{fig:Whitney}.

To construct the structures of Theorem~\ref{thm:defstr} with
parameters $(s,t)\in U$, we need to construct a family of
representations $\rho_{(s,t)}$ of $\pi_1M$ in
$\operatorname{Isom}^+(\mathbb X^3)$, which are going to be the
holonomy representations of the structures. In fact we work in the
universal covering
 of $\operatorname{Isom}^+(\mathbb X^3)$, that  we denote by $G$.
 When $\mathbb X^3=\mathbb H^3$ then  $G=SL_2(\mathbb C)$, and
when $\mathbb X^3=\mathbb S^3$ then
 $G=SU(2)\times SU(2)$.

The starting point in the construction of $\rho_{(s,t)}$ is the
holonomy representation
\[
\hol\co\pi_1M\to \operatorname{Isom}(Nil)
\]
and the exact sequences:
\[
    0\to \mathbb R\to
\operatorname{Isom}(Nil)
    \mathop{\longrightarrow}^{\pi} \operatorname{Isom}(\R^2)\to 1,
  \]\[
     0\to \mathbb R^2\to \operatorname{Isom}(\R^2)
\mathop{\longrightarrow}^{\ROT}  O(2)\to 1.
\]
The first one comes from the Riemannian fibration $\mathbb R\to
Nil\to\mathbb R^2$ and the second one is well known. We consider the
representation
 \[
    \phi_0=\ROT\circ\pi\circ\hol\co\pi_1M\to O(2)\subset SO(3)
 \]
and we lift it to $\rho_0\co\pi_1M\to SU(2)\cong
\widetilde{SO(3)}$. We fix $x_0\in\mathbb X^3$ and we view $SU(2)$
as the stabilizer of $x_0$ in $G$. We construct $\rho_{(s,t)}$ as
a perturbation of $\rho_0$. The infinitesimal properties of this
perturbation are related to the holonomy representation $\hol$ and
to sections to the above exact sequences, because by composing
$\hol$ with those sections we obtain cocycles and cochains.

\paragraph{Organization of the paper}
We start with a review of $Nil$ geometry and the holonomy
representation in Section~\ref{sect:nil}, pointing out its
cohomological aspects for relating it later to infinitesimal
deformations. In Section~\ref{sec:def} we construct the
deformation spaces for spherical and hyperbolic structures,
assuming the existence of suitable representations $\rho_{(s,t)}$.
Those representations are constructed in Section~\ref{sec:varrep},
and their infinitesimal properties are studied in
Section~\ref{sect:inf}. Section~\ref{sec:euclidean} is devoted to
Euclidean structures, obtained as degeneration of hyperbolic and
spherical ones. In Section~\ref{sect:DF} we analyze the Dehn
filling parameters, achieving the proof of Theorem~B. The part of
Theorem~A not contained in Theorem~B is proved in
Section~\ref{sect:path}, together with
 Propositions~\ref{prop:properties1}
and~\ref{prop:properties2}. Section~\ref{sec:example} is devoted to
an example, where the limit $l_0$ of
Propositions~\ref{prop:properties1} and~\ref{prop:properties2} is
explicitly computed. Finally Section~\ref{sect:cohmology} is devoted
to the proof of some technical computations in cohomology.

\section{The holonomy representation}~\label{sect:nil}
The usual model for $Nil$
 is the Heisenberg group of matrices of the
form
\[
Nil=\left\{ \left(\left.
    \begin{array}{ccc}
    1 & x & z \\
    0 & 1 & y \\
    0 & 0 & 1
    \end{array}
\right)\right\vert\, x,y,z\in\R  \right\}
\]
which is canonically identified to $\R^3$ by taking coordinates
$(x,y,z)$. For our purposes it will be convenient to work with
another model. Following \cite{Sua}, we consider $\R^3$ with the
product:
\[
(x_1,x_2,x_3)(y_1,y_2,y_3)=(x_1+y_1, x_2+y_2,
 x_3+y_3+ x_1 y_2- x_2 y_1)
\]
This is another model for $Nil$. The isomorphism between both models
is given by $x=\sqrt2 x_1$, $y=\sqrt2 x_2$ and $z=x_3+x_1 x_2$.

\subsection{The isometry group of $Nil$}

We consider a $2$--parameter family of left-invariant metrics
\[ ds^2=\lambda^2 (d x_1^2+d x_2^2)+ \mu^2 (d x_3+x_2\,d  x_1 -
x_1\, d x_2)^2
\]
 for $\lambda,\mu\in\R-\{0\}$.
 All these metrics have the
same 4--dimensional isometry group $\Isom(Nil)$. This group
$\Isom(Nil)$ preserves the orientation, it has two components, and it
is a semi-direct product
\[
\operatorname{Isom}(Nil) \cong Nil\rtimes O(2).
\]
The group $O(2)$ acts on $Nil$ linearly  as the projection of the
standard action of $O(2)\subset SO(3)$  on $\R^3\cong Nil$ preserving
the plane $x_3=0$. To see that this action is an isometry, it may be
useful to write the metric in cylindrical coordinates
$x_1=r\cos\theta$ and $x_2=r\sin\theta$:
\[
ds^2=\lambda^2(d r^2+r^2 d\theta^2)+\mu^2(d x_3-r^2 d\theta)^2.
\]
The projection $Nil\to\R^2$ that maps $(x_1,x_2,x_3)\in Nil$ to
$(x_1,x_2)\in\R^2$ is a Riemannian fibration with fibre a line $\R$.
This fibration is preserved by the isometry group and induces an
exact sequence
 \begin{equation}\label{eqn:vert}
    0\to \mathbb R\to \operatorname{Isom}(Nil)
    \mathop{\longrightarrow}^{\pi} \operatorname{Isom}(\R^2)\to 1.
\end{equation}
A section
    $$ \VERT_p\co \Isom(Nil)\to \R $$
may be constructed by fixing a base point $p\in Nil$ as follows:
for any $g\in \Isom(Nil)$, $\VERT_p(g)$ is the third coordinate of
$g(p)p^{-1}$.

 On the other hand, we have the well known split exact sequence
 \begin{equation}\label{eqn:trans}
0\to \mathbb R^2\to \operatorname{Isom}(\R^2)
\mathop{\longrightarrow}^{\ROT}  O(2)\to 1.
 \end{equation}
 A section
 \[
 \TRANS_q\co \Isom(\R^2)\to\R^2
 \]
  may also be constructed fixing
 a base point $q\in\mathbb R^2$. For any $q\in \Isom (\R^2)$,
 $\TRANS_q(g)=g(q)-q\in \mathbb R^2$.

\subsection{The holonomy representation}\label{subs:holrep}

 Our starting point is the  holonomy
representation of the orbifold $\OO$: \[\hol\co \pi_1^o(\mathcal
O)\to \Isom(Nil)\] and the representation induced on the open
manifold $M=\vert \OO\vert-\Sigma$.

\begin{defn} Given the induced
representation $\hol\co\pi_1M\to\Isom(Nil)$, $p\in Nil$ and
$q=\pi(p)\in\mathbb R^2$,
 we
define the following maps:
\begin{eqnarray*}
 \phi_0= \ROT\circ \pi\circ \hol\co \pi_1M &\to& O(2)\subset SO(3),\\
 z_q= \TRANS_q\circ\pi\circ\hol\co\pi_1M &\to& \mathbb R^2,\\
c_p=\VERT_p\circ\hol\co \pi_1M &\to&\mathbb R.
\end{eqnarray*}
\end{defn}

Those three maps determine uniquely the representation $\hol$. It is
clear that $\phi_0$ is also a representation, but $z_q$ and $c_p$ are
not. However  they satisfy some cohomological conditions that we
describe next.
 To do it, we view
both $\mathbb R^2=\mathbb R^2\times 0$ and $\mathbb R=0\times\mathbb
R$ as subspaces of $\mathbb R^3$, therefore they are $\pi_1M$--modules
via $\phi_0 \co \pi_1M\to O(2)\subset SO(3)$.

The map $z_q$
 is a cocycle twisted by $\phi_0$. This is,
    \begin{equation}\label{eqn:cocycle}
z_q(g_1g_2)=z_q(g_1)+\phi_0(g_1)z_q(g_2),\qquad \forall
g_1,g_2\in\pi_1M
    \end{equation}
The map $c_p$ satisfies the following relation:
\[
c_p(g_1g_2)-c_p(g_1)-\phi_0(g_1) c_p(g_2)= z_q(g_1)\times
\phi_0(g_1)z_q(g_2),\qquad \forall g_1,g_2\in\pi_1M,
\]
 where $\times $ denotes the usual cross
product in $\mathbb R^3$. In cohomology terms, the previous
inequality is:
\[
\delta(c_p)=z_q\cup z_q,
\] where $\delta$  denotes the cobundary, and $\cup$, the cup product
associated to $\times$.

The set of all cochains (ie, maps $\pi_1M\to \mathbb R^2\times 0$)
is a vector space denoted by $C^1(\pi_1M,\mathbb R^2\times 0)$. The
subspace of all cocycles (ie, maps $\pi_1M\to \mathbb R^2\times 0$
satisfying (\ref{eqn:cocycle})) is  denoted by $Z^1(\pi_1M,\mathbb
R^2\times 0)$. Hence $z_q\in Z^1(\pi_1M,\mathbb R^2\times 0)$.

Let $B^1(\pi_1M,\mathbb R^2\times 0)$ denote the subspace of all
coboundaries, ie, cocycles $b_r$ with the property that there exists
$r\in \mathbb R^2\times 0$ with $b_r(g)=r-\phi_0(g)(r)$, $\forall
g\in\pi_1M$. The cocycle $ z_q\not\in B^1(\pi_1M,\mathbb R^2\times 0)
$ because $z_q$ does not have a global fixed point in $\mathbb
R^2\times 0$. Thus the cohomology class of $z_q$ in
\[
H^1(\pi_1M,\mathbb R^2\times 0)=Z^1(\pi_1M,\mathbb R^2\times
0)/B^1(\pi_1M,\mathbb R^2\times 0)
\]
is not zero, and it may be easily checked that it is independent of
the choice of $q\in \mathbb R^2\times 0$.

 We will prove at the end of
the paper that
\[
H^1(\pi_1M,\mathbb R^2\times 0)\cong\mathbb R \qquad\textrm{ and
}\qquad H^1(\pi_1M,0\times \mathbb R)\cong 0.
\]
This has two consequences. Firstly $z_q$ is unique up to the choice
of $q$ and up to homoteties. Secondly, once $z_q$ and
$p\in\pi^{-1}(q)$ have been fixed,   then $c_p$ is unique.

Different choices of the cohomology class $[z_p]$ correspond to the
composition of the holonomy with an automorphism:
\[
\begin{array}{rcl}
    Nil&\to&Nil\\
    (x_1,x_2,x_3)&\mapsto & (\lambda x_1,\lambda x_2,\lambda^2 x_3)
\end{array}
\]
 for some $\lambda\in\mathbb R-\{0\}$.


\subsection{Lifting the holonomy}

 We recall that $M=\vert O\vert-\Sigma$ and that
$
 \phi_0=\ROT\circ \pi\circ \hol\co \pi_1 \mathcal O\to O(2)\subset SO(3)
 $. The
 representation  of $\pi_1M$ in $SO(3)$ induced by
$\phi_0$ lifts to a representation to
 $SU(2)\cong Spin(3)=\widetilde{SO(3)}$,
 because we can view it as the holonomy of a non-complete structure
 on $M$ and apply the following result of Culler \cite{Cul}.

\begin{lem}{\rm\cite{Cul}}\label{lem:cul}\qua A spin structure on 
$M$ determines a lift of $ \ROT\circ \pi\circ \hol$ to $Spin(3)\cong
SU(2)$. In particular, since $dim(M)=3$ there exists a lift.\qed
\end{lem}

\begin{rem}
Two spin structures determine a morphism $\theta\co \pi_1M\to
\Z/2\Z$. It follows from the construction of \cite{Cul}, that if
$\rho_1$ and $\rho_2$ are the lifts associated to these structures,
then
\[
\rho_1(g)= (-1)^{\theta(g)}\rho_2(g)\qquad\textrm{ for every
}g\in\pi_1M.
\]
\end{rem}

 From now on we fix a spin structure on $M$, hence we also fix a lift
of $ \phi_0=\ROT\circ \pi\circ \hol$:
 \[
        \rho_0\co \pi_1M\to SU(2).
 \]

\subsection{Changing the spin structure}

We consider the natural surjection
\[
\theta\co \pi_1M\twoheadrightarrow \Z/2\Z
\]
which is the composition of $\phi_0\co \pi_1M\to O(2)$ with the
projection $O(2)\twoheadrightarrow\pi_0(O(2))\cong \Z/2\Z $.

We consider the change of spin structure associate to $\theta$.

If $\rho$ is the
 lift of a representations of $\pi_1M$ in
$\operatorname{Isom}^+(\mathbb X^3)$ as in 
Lemma~\ref{lem:cul}, this change of spin structure corresponds to
to replace the lift $\rho$ by $(-1)^\theta\rho$.

\begin{lem}\label{lem:spinrho0}
The representation $(-1)^\theta\rho_0$ is conjugate to $\rho_0$.
\end{lem}

\proof It suffices to check that
$\operatorname{trace}((-1)^\theta\rho_0(g))=\operatorname{trace}(\rho_0(g))$,
for every $g\in\pi_1M$, because $\rho_0$ is a representation in
$SU(2)$.  If $g\in\ker\theta$, then the equality of traces holds true
because $(-1)^{\theta(g)}\rho_0(g)=\rho_0(g)$. If $\theta(g)=1$ then
$\ROT\circ\pi\circ\hol(g)$ is a rotation of angle $\pi$, as every
element in $O(2)-SO(2)$ viewed in $SO(3)$. Hence
$\operatorname{trace}(\rho_0(g))=0$ and therefore
$\operatorname{trace}((-1)^\theta\rho_0(g))=-\operatorname{trace}(\rho_0(g))=0$.
 \qed

\section{Deformation spaces}\label{sec:def}

From now on $\mathbb X^3$ will denote $\mathbb H^3$ and $\mathbb
S^3$. Every statement about $\mathbb X^3$ will be understood to be a
statement about both, the hyperbolic space and the $3$--sphere. The
hyperbolic plane and the $2$--sphere will be denoted by $\mathbb X^2$.

\subsection{Spaces of geometric  structures}

\begin{thm}\label{thm:defstr}
 There exists a space of geometric structures on $M=\OO-\Sigma$ with
 Dehn filling end
 parametrized by a neighborhood of the origin $U\subset \mathbb R^2$.
 According to the parameters $(s,t)\in U$, the structure is
 of the following kind:
\begin{enumerate}[\rm(i)]
    \item the original $Nil$ structure, when $(s,t)=0$;
    \item modeled on $\mathbb X^3$, when $s\, t\neq 0$;
    \item a foliation transversely modelled on $\mathbb X^2$, when
    $s=0$, $t\neq 0$; and
    \item a Euclidean structure, when $s\neq 0$, $t=0$.
\end{enumerate}
In addition, those structures are oriented and equipped with a spin
structure, so that $(s,- t)$ and $(s,t)$ correspond to structures
with opposite orientation, and $-(s,t)$ and $(s,t)$ correspond to
 the spin structures differing by $\theta$.
\end{thm}

A Dehn filling end for the  structure on  $T^2\times (0,1]$ means the
following. There is a geodesic $\gamma\subset\mathbb X^3$ such that
the developing map $D\co \widetilde{T^2}\times  (0,1]\to\mathbb X^3$
maps $\{x\}\times (0,1]$ to a minimizing segment between $D(x,1)$ and
$\gamma$, for every $x\in\widetilde{T^2}$. In addition, the parameter
in $ (0,1]$ is proportional to arc-length.

The parameter $(s,t)$ has the following interpretation. We choose
 $l, m\in\pi_1M$ so that they generate a peripheral group and $m$
 is a meridian for $\Sigma$. We may choose $l$ so that $\phi_0(l)$ is
 trivial.
The rotation angle and the translation length of the
 holonomy of $l$ are respectively $s$ and $t$.

\begin{conv}
We fix $x_0\in\mathbb X^3$ and we view $\phi_0$ as a representation
in $\Isom^+(\mathbb X^3)$ that fixes $x_0$, because  $SO(3)$ is the
stabilizer of a point in $\Isom^+(\mathbb X^3)$. We also fix
$\{e_1,e_2,e_3\}$ a positive orthonormal basis for $\mathbb R^3$ so
that $\langle e_1,e_2\rangle=\mathbb R^2\times 0$ and  $\langle
e_3\rangle=0\times \mathbb R$ are the subspaces invariant by $O(2)$.
The totally geodesic plane tangent to $\mathbb R^2\times 0\subset
T_{x_0}\mathbb X^3$ is denoted by $\mathbb X^2=\exp_{x_0}(\mathbb
R^2\times 0)$.
\end{conv}

The following maps from $Nil$ to $\mathbb X^3$ will be used in the
proof of Theorem~\ref{thm:defstr}.

\begin{defn} For $(s,t)\in\R^2$ we define:
\[\begin{array}{rcl}
  \Delta_{(s,t)}\co Nil\cong\R^3 & \to & \mathbb X^3 \\
  (x_1,x_2,x_3) & \mapsto & \exp_{x_0}(t(x_1 e_1+x_2e_2+s\ x_3e_3))
\end{array}
\]
where $\exp_{x_0}$ denotes the Riemannian exponential at the point
$x_0\in\mathbb X^3 $. Here we have identified $Nil$ with $\mathbb
R^3$.
\end{defn}

Notice that, when $s\, t\neq 0$, $\Delta_{(s,t)}$ is a local
diffeomorphism, and when $s=0$ but $t\neq 0$, it is a local a
submersion of rank 2 onto $\mathbb X^2$.

\subsection{Deformations of representations}

\begin{prop}\label{prop:defrep}
 There exists a perturbation $\rho_{(s,t)}\co\pi_1M\to
G$ of $\rho_0$, with parameter $(s,t)\in U\subset \mathbb R^2$, such
that:
\begin{enumerate}[\rm(i)]
\item $\rho_{(s,0)}$ stabilizes $x_0$.
\item $\rho_{(0,t)}$ stabilizes $\mathbb X^2=\exp_{x_0}(\mathbb R^2\times 0)$
\item For every $g\in\pi_1M$
\[
\lim_{\begin{smallmatrix}
  (s,t)\to 0 \\
  s\, t\neq 0
\end{smallmatrix}} \Delta_{(s,t)}^{-1} \circ \rho_{(s,t)}(g)\circ \Delta_{(s,t)} =
\hol(g)
\]
uniformly on compact subsets of $Nil$ for the $\mathcal
C^1$--topology.
\item Let $\overline\Delta_t=\Delta_{(0,t)}\vert_{\mathbb R^2\times 0}$.
For every $g\in\pi_1M$
\[
\lim_{t\to 0} \overline \Delta_{t}^{-1} \circ \rho_{(0,t)}(g)\circ
\overline \Delta_{(0,t)} =\pi\circ \hol(g)
\]
uniformly on compact subsets of $\mathbb R^2\times 0$ for the
$\mathcal C^1$--topology.
\item The representations $\rho_{(-s,-t)}$ and $(-1)^{\theta}\rho_{(s,t)}$ are
conjugate in $\widetilde{\operatorname{Isom}^+}(\mathbb X^3)$.
\item  $\rho_{(-s,t)}$ and $\rho_{(s,t)}$ are
conjugate by an orientation reversing element in
$\widetilde{\operatorname{Isom}}(\mathbb X^3)$.
\end{enumerate}
\end{prop}

We shall prove Theorem~\ref{thm:defstr} assuming this proposition.
The perturbation we will construct satisfy some more properties
related to the Euclidean structures, when $t=0$. These properties
will be explained later, hence for the moment  we will not prove the
part of Theorem~\ref{thm:defstr} concerning Euclidean structures.

Properties (iii) and (iv) of the proposition are related to the
infinitesimal properties of $\rho_{(s,t)}$ and to the cocycle $z_q$
and the cochain $c_q$.

\subsection{Proof of Theorem~\ref{thm:defstr}}

We construct a covering $\{U_i\}_{i=0,\ldots,n}$ of $M$  such that
$U_i$ is 1--connected for $i\geq 1$ and $U_0$ is a neighborhood of the
end of $M$.

Since $U_1$ is simply connected, the lift of $U_1$ in the universal
covering of $M$ is
\[
\widetilde{U_1}=\bigcup_{g\in\pi_1M} g W_1
\]
for some open set $W_1\subset \widetilde M$ that projects
homeomorphically to $U_1$. We define on $W_1$
\[
D_{(s,t)}\vert_{W_1}=\Delta_{(s,t)}\circ \mathcal D_0\vert_{W_1}\co
W_1\to \mathbb X^3
\]
Where $\mathcal D_0\co \widetilde M\to Nil$ is the holonomy for
the Nil structure. Next we define $D_{(s,t)}$ on $\widetilde{U_1}$
by taking the equivariant extension. By
Proposition~\ref{prop:defrep} (iii),
\[
\lim_{\begin{smallmatrix}
  (s,t)\to 0 \\
  s\, t\neq 0
\end{smallmatrix}} \Delta_{(s,t)}^{-1} \circ D_{(s,t)}\vert_{\widetilde{U_1}} =
\mathcal D_0\vert_{\widetilde{U_1}}
\]
for the $\mathcal C^1$--topology uniformly on compact subsets.

We make the same construction for all sets $U_i$ with $i\geq 2$ and
for $U_0$ we make a cylindrical construction taking care of the
holonomy at the end. We glue the final construction by using standard
techniques about bump functions and refinements, as explained in
\cite{CEG} (and also in \cite{PoTwo,HPS}), so we obtain a
 family $D_{(s,t)}$ of maps that are $\rho_{(s,t)}$--equivariant and
 such that:
 \[
\lim_{\begin{smallmatrix}
  (s,t)\to 0 \\
  s\, t\neq 0
\end{smallmatrix}} \Delta_{(s,t)}^{-1} \circ D_{(s,t)} =
\mathcal D_0
\]
for the $\mathcal C^1$--topology  uniformly on compact subsets. In
particular $D_{(s,t)}$ is a local diffeo for small values of $(s,t)$
with $s\, t\neq 0$.

For the neighborhood $U_0$ we need to be more careful. Let
$\gamma=\exp_{x_0}\langle e_1\rangle$ be the geodesic preserved by
$\rho_0(m)$. We will construct $\rho_{(s,t)}$ so that $\gamma$ is
preserved by $\rho_{(s,t)}(m)$ (and also by $\rho_{(s,t)}(l)$, by
commutativity). It will also follow from the construction that
$\rho_{(s,t)}(l)$ is the composition of a translation of length
$t$ with a rotation of angle $s$ around $\gamma$. We consider the
family of maps $\Theta_{(s,t)}\co \tilde U_0\to \mathbb
R^3-\langle e_1\rangle$ such that $\Theta_0$ is the developing map
$\mathcal D_0$ restricted to $\tilde U_0$, the distance from
$\Theta_{(s,t)}(x)$ to $\langle e_1\rangle$ is independent of
$(s,t)$ and $\Theta$ is $\pi_1U_0$--equivariant by the action of
$\rho_{(s,t)}$. Then we define $\mathcal D_{(s,t)}\vert_{\tilde
U_0}= \Delta_{(s,t)}\circ \Theta_{(s,t)}$ and we glue it in the
same way.

This proves assertion (ii) of Theorem~\ref{thm:defstr}. The proof
of assertion (iii) is quite similar by using
Proposition~\ref{prop:defrep}. The properties about symmetries are
also clear from Proposition~\ref{prop:defrep}.
 \endproof

We recall that the part of the theorem concerning Euclidean
structures will be proved later.

\section{Construction of the representations}\label{sec:varrep}

In this section we construct the representations of
Proposition~\ref{prop:defrep}.

\subsection{Smoothness of the varieties of characters}
 We work
with the varieties of representations of $\pi_1M$ in $SU(2)$ and
$SL_2(\mathbb C)$:
\begin{eqnarray*}
R(M,SU(2)) &=&\Hom(\pi_1M, SU(2)),
\\
 R(M,SL_2(\C)) &=&\Hom(\pi_1M, SL_2(\C)).
\end{eqnarray*}
The varieties of characters are defined as:
\begin{eqnarray*}
    X(M,SU(2)) &=&R(M,SU(2)) /SU(2),
    \\
    X(M,SL_2(\C)) &=&R(M, SL_2(\C))/\!/ SL_2(\C).
\end{eqnarray*}
The symbol $/\!/$ in the definition of $X(M,SL_2(\C))$ means the
algebraic quotient (in invariant theory). In particular
$X(M,SL_2(\C))$ is algebraic affine (also defined over $\Q$). However
since $SU(2)$ is compact but not complex, $X(M,SU(2))$ is just the
topological quotient, and it is only real semi-algebraic, contained
in the set of real points of $X(M,SL_2(\C))$.

Every point in $X(M,SL_2(\C))$ is the character of a representation
in $SL_2(\C)$, ie, a map
\[
    \begin{array}{rcl}
    \chi_{\rho}\co \pi_1M&\to& \C\\
    \gamma&\mapsto&\operatorname{trace}({\rho}(\gamma))
    \end{array}
\]
 for some
$\rho\in R(M, SL_2(\C))$. Every conjugacy class of representation
into $SU(2)$ is determined by its character, therefore the notation
makes sense and $X(M,SU(2))\subset X(M,SL_2(\C))$.

\begin{defn} For every $\gamma\in\pi_1M$,
$I_{\gamma}\co X(M,SL_2(\C))\to\C$ denotes the evaluation map. In
other words, it is the map induced by the trace function:
\[
I_{\gamma}(\chi_{\rho})=\chi_{\rho}(\gamma)=\operatorname{trace}(\rho(\gamma)).
\]
\end{defn}

\begin{prop}\label{prop:smooth} The character $\chi_0$ of $\rho_0$ is a smooth
one dimensional point  of both $X(M,SU(2))$ and $X(M,SL_2(\mathbb
C))$.
\end{prop}

\proof We first prove the proposition for $X(M,SL_2(\mathbb C))$.
 By a Theorem~5.6 of Thurston's notes \cite{ThurstonNotes},
the local dimension of $X(M,SL_2(\mathbb C))$ at the character of
$\rho_0$ is at least one. It suffices to prove that
$H^1(\pi_1M,sl_2(\mathbb C))\cong \mathbb C$, (where $\pi_1M$ acts on
$sl_2(\mathbb C)$ via $Ad\rho_0$) because this cohomology group
contains the Zariski tangent space of $X(M,SL_2(\mathbb C))$ at
$\chi_0$. We have said before that $H^1(\pi_1M,\mathbb R^2\times
0)\cong \mathbb R$ and $H^1(\pi_1M,0\times\mathbb R)\cong 0$.
Therefore
\[
H^1(\pi_1M,su(2))\cong H^1(\pi_1M,\mathbb R^3)\cong\mathbb R,
\]
because $su(2)$ and $\mathbb R^3$ are isomorphic as $\pi_1M$--modules.
In particular \[ H^1(\pi_1M,sl_2(\mathbb C))\cong
H^1(\pi_1M,su(2))\otimes_{\mathbb R}\mathbb C\cong\mathbb C.
\]
The proposition for $X(M,SU(2))$ follows easily, using the fact that
the variety $X(M,SL_2(\mathbb C))$ is defined over $\mathbb R$ and a
neighborhood of $\chi_0$ in $X(M,SL_2(\mathbb C))\cap\mathbb R^N$
coincides with $X(M,SU(2))$.\qed

\subsection{Local parametrization}\label{subsec:localpara}

We construct a local parameter of a neighborhood of $\chi_0$ in
$X(M,SL_2(\C))$. We choose $l, m\in\pi_1M$ so that they generate a
peripheral subgroup $\pi_1T^2$. We assume that
 $m$ is a meridian of
$\Sigma$. We also assume that $\theta(l)=0$, by replacing $l$ by
$l\,m$ if necessary.

\begin{rem} We have that $\rho_0(l)=\pm Id$, because $l$ and $m$
commute,
and $l\in\ker\theta$ but $\theta(m)=1$ (ie, $\phi_0(l)\in SO(2)$ but
$\phi_0(m)\in O(2)-SO(2)$ ).
\end{rem}

The idea is to choose $w=\alpha_l$ the angle rotation of $\rho(l)$
as a local parameter of $X(M,SU(2))$ (so that its extension to
$X(M,SL_2(\C))$ corresponds to $\pm 2i$ times the logarithm of an
eigenvalue). The sign of this angle is determined by the sense of
rotation around the invariant geodesic, which corresponds to a
choice of the spin structure and determines the choice of the lift
$\rho_0$. We would like to define $w$ as $2\arccos (I_l/2)$, but
$\arccos$ is not well defined in a neighborhood of $\pm 1$.
Formally, we can define it as follows.

\begin{defn} \label{def:w} In a neighborhood of $\chi_0$ we define $w$ as
\[
w=2\arccos( I_{l m}/2)-2\arccos( I_m/2).
\]
so that $I_l= \pm 2\cos \frac{w}2$.
\end{defn}

\begin{lem}\label{lem:w}
 The function $w$ defines a local parametrization of both
varieties of characters $X(M,SU(2))$ and $X(M,SL_2(\C))$.
\end{lem}

\proof
 It follows from the proof of Proposition~\ref{prop:smooth}
that $H^1(\pi_1M,su(2))$ is isomorphic to the tangent space
$T_{\chi_0}X(M,SU(2))$. Thus we view $H_1(\pi_1M,$ $su(2))$ as the
cotangent space $T_{\chi_0}^1X(M,SU(2))\cong \mathbb R$, and it is
sufficient to check that the differential form $d\, w\neq 0$. In
particular, we just need to prove that the Kronecker pairing $\langle
d\, w,z_q\rangle$ does not vanish, where $z_q$ is the cocycle defined
in Subsection~\ref{subs:holrep}. Since, for a representation $\rho$,
$w(\chi_{\rho})$ is precisely the angle of $\rho(l)$, Proposition~9.6
in \cite{PoOne} implies that $\langle d\, w,z_q\rangle$ is precisely
the translation length of $\pi\circ\hol(l)$. This length is non-zero
because $\Sigma$ is horizontal.\endproof

We recall that
 $ \theta\co \pi_1M\twoheadrightarrow \Z/2\Z
 $
 is the composition of $\phi_0\co \pi_1M\to O(2)$ with the
projection $O(2)\twoheadrightarrow\pi_0(O(2))\cong \Z/2\Z $. We
consider the change of spin structure associate to $\theta$. For a
representation $\rho\in R(M,SL_2(\C))$, to change the spin structure
corresponds to replace $\rho$ by $(-1)^\theta\rho$.

\begin{lem}\label{lem:spin}
 $w(\chi_{(-1)^\theta\rho})=-w(\chi_{\rho})$.
\end{lem}

\proof Since $\chi_0$ is invariant by $\theta$
(Lemma~\ref{lem:spinrho0}),
 the neighborhood of
$\chi_{0}$ 
may be chosen invariant by the
change of the spin structure. Since $\theta(m)=\theta(lm)=1$, we
have that $I_m(\chi_{(-1)^\theta\rho})=- I_m(\chi_{\rho})$ and
$I_{lm}(\chi_{(-1)^\theta\rho})=- I_{lm}(\chi_{\rho})$. Therefore,
for the branch of $\arccos$ with $\arccos(0)=\pi/2$, we have:
\[\begin{array}{rcl}
  2\arccos(-I_m(\chi)/2)    & = & \pi-2\arccos(I_m(\chi)/2) \\
  2\arccos(-I_{lm}(\chi)/2) & = & \pi-2\arccos(I_{lm}(\chi)/2)
\end{array}
\]
and the lemma follows.\qed

\subsection{Deformations of characters}\label{subsec:defchar}

We choose different varieties of  characters for the hyperbolic and
the spherical case, but we will unify the notation for the
neighborhood $U$.

In the hyperbolic case, since $\Isom^+(\H^3)\cong PSL_2(\C)$  we work
in $X(M,SL_2(\C))$ (we recall that we have fixed a spin structure,
hence all holonomy representations have a natural lift). We fix
$U\subset \R^2$ a neighborhood of the origin, with coordinates
$(s,t)\in U$ and set
    \[
    w=s-t\, i,
    \]
so that, for any representation $\varrho_w$  with character $w$, the
complex length of $\varrho_w(l)$ is  $i\, w=t+s\, i$ (ie, a
translation of length $t$ plus a rotation of angle $s$).

 In the
spherical case, since $Spin(4)\cong SU(2)\times SU(2)$ we work in
$X(M,SU(2))\times X(M,SU(2))$. We denote by $w_1$ and $w_2$ the
ordered (real) parameteters of each factor $X(M,SU(2))$ given by
 Definition~\ref{def:w}.
We fix $U\subset \R^2$ a neighborhood of the origin, with coordinates
$(s,t)\in U$ and we set
    \[ (w_1,w_2)=(s+t,s-t)\]
Again, any representation   with character $(w_1,w_2)$ evaluated at
$l$ is a translation of length $t$ composed with a rotation of angle
$s$ around the same edge.

In both cases, $\chi_0$ the character of $\rho_0$ has coordinates
$(s,t)=(0,0)$. To construct the representations $\rho_{(s,t)}$ we
need a section to the projection  \[R(M,SL_2(\mathbb C))\to
X(M,SL_2(\mathbb C)).\] This will be done after the description of
$\rho_0$.

\subsection{Description of $\rho_0$}

We recall that we have fixed $\{e_1,e_2,e_3\}$ an orthonormal basis
for $\mathbb R^3$, so that $\langle e_1, e_2\rangle=\mathbb R^2\times
0$ and $\langle e_3\rangle=0\times \mathbb R$ are the subspaces
invariant by $\rho_0$.

By using the natural identification $su(2)\cong \mathbb R^3$ as
$SU(2)$--modules, we view $e_1,e_2,e_3$ as three matrices of $su(2)$
such that the following formula and its cyclic permutations hold:
\[
[e_1,e_2]=e_3,
\] because the Lie
bracket in $su(2)$ corresponds to the cross product in $\R^3$.

\begin{rem} Let $v\in su(2)\cong \R^3$ be a unitary vector and
$\alpha\in\R$. Then $\exp(\alpha v)\in SU(2)$ projects in $SO(3)$ to
a rotation of angle $\alpha$ around $\langle v\rangle$.
\end{rem}

Thus, if $\theta(g)=0$ then $\rho_0(g)=\exp(\alpha_g e_3)$, for some
$\alpha_g\in\mathbb R$.
 Notice
 that
    \[
    \exp((\alpha_g+2\pi)e_3)=-\exp(\alpha_g e_3) .
    \]
  We may also
assume that $e_1$ is the vector invariant  for the meridian $m$ and
that the spin structure has been chosen so that $\rho_0(m)=\exp(\pi
e_1)$. The elements which are not in the kernel of $\theta$ are of
the form $g m$ for some $g\in\ker(\theta)$, and we have $\rho_0(g m
)=\exp(\pi(\cos(\alpha_g/2) e_1+ \sin(\alpha_g/2) e_2))$.

\begin{rem}\label{rem:rho0inv}  The
conjugation matrix between $\rho_0$ and $(-1)^{\theta}\rho_0$ is
$\pm\exp(\pi e_3)$.
\end{rem}

This remark follows from the description of $\rho_0$ and the fact the
adjoint action (equivalent to the orthogonal action on $\R^3$) of
$\exp(\pi e_3)$ changes the sign of $e_1$ and $e_2$ and preserves
$e_3$.

\subsection{The section for $R(M,SL_2(\C))$}

\begin{lem}\label{lem:sigma0} There exists a neigborhood $V\subset X(M,SL_2(\C))$ and a
section $\sigma\co V\to R(M,SL_2(\C))$ such that,
    if $\varrho_{w}=\sigma(w)$, then $\forall g\in\pi_1M$,
    \[  \varrho_{w}(g)=
    \exp(f_{g}(w)+h_g(w))\rho_0(g)
    \]
where $f_g$ and $h_g$  are analytic maps with real coefficients
valued on the Lie algebra $sl_2(\C)$, such that $f_g(w)\in \langle
e_1,e_2\rangle_{\C}$, $h_g(w)\in \langle e_3\rangle_{\C}$, $f_g$ is
odd and $h_g$ is even.
\end{lem}

When we say that the coefficients  of $f_g$ and $h_g$ are real, we
mean that for $w\in\mathbb R$, $f_g(w),h_g(w)\in su(2)$.

\proof The proof is based in a construction analogue to Luna's slice
theorem. We consider the involution $\nu$ on $R(M,SL_2(\C))$ and
$R(M,SU(2))$ defined as follows:
\[
\nu(\rho)=(-1)^{\theta}Ad_{\exp(\pi e_3)}\circ\rho
\]
where $\theta\co \pi_1M\twoheadrightarrow\mathbb Z/2\mathbb Z$ is
described above.

 By the remark in Subsection~\ref{rem:rho0inv}, $\nu(\rho_0)=\rho_0$. In addition, by
Lemma~\ref{lem:spin}, if $t\co R(M,SL_2(\C))\to X(M, SL_2(\C))$
denotes the projection, then
\[
    w\circ t\circ\nu=-w\circ t.
\]

\begin{lem}\label{lem:S} There exists an algebraic complex curve
$\mathcal S\subset R(M,SL_2(\C))$ with the following properties:
    \begin{enumerate}[\rm(i)]
    \item $\rho_0$ is a smooth point of $\mathcal S$.
    \item The projection $t\co R(M,SL_2(\C))\to X(M,SL_2(\C))$ restricts
    to a map $t\vert_{\mathcal S}\co \mathcal S\to X(M,SL_2(\C))$
    locally bianalytic at $\rho_0$.
    \item $\mathcal S'=\mathcal S\cap R(M,SU(2))$ is a real curve smooth at
    $\rho_0$ and the restriction $t\vert_{\mathcal S'}\co \mathcal S'\to X(M,SU(2))$
    is also locally bianalytic at $\rho_0$.
    \item $\mathcal S$ is invariant by the involution $\nu$.
    \item For every $\rho\in\mathcal S$,
    $\rho(m)=\exp(\alpha_{\rho}\,
    e_1)$, for some $\alpha_{\rho}\in\mathbb R$.
    \end{enumerate}
\end{lem}

We postpone its proof. Assuming it holds, we conclude the proof of
Lemma~\ref{lem:sigma0}. It suffices to take $\sigma=t\vert_{\mathcal
S}^{-1}$.  We write $\varrho_{w}=\sigma(w)$ and $\varrho_{w}(g)=
    \exp(f_{g}(w)+h_g(w))\rho_0(g)
$
for some analytic maps such that the image of $f_g$ is contained in
$\langle e_1,e_2\rangle_{\C}$ and the image of $h_g$ is contained in
$\langle e_3\rangle_{\C}$. These maps have real coefficients by
assertion  (iii) of Lemma~\ref{lem:S}. We use the involution to prove
that $f_g$ is odd and $h_g$ is even. The representations
$\varrho_{-w}$ and $\nu(\varrho_{w})$ have the same character. By the
properties of $\mathcal S$, it follows that
$\varrho_{-w}=\nu(\varrho_{w})$. In addition
    {\setlength\arraycolsep{2pt}
\begin{eqnarray}
    \nu(\varrho_{w})(g)  &=& (-1)^{\theta(g)}Ad_{\exp(\pi
    e_3)}(\varrho_{w}(g))\nonumber \\
                    &=&(-1)^{\theta(g)}Ad_{\exp(\pi
    e_3)}(\exp(f_{g}(w)+h_g(w)))\ Ad_{\exp(\pi e_3)}(\rho_0(g))
    \nonumber \\
    \label{eqn:nu} &  = & \exp(-f_g(w)+h_g(w))\rho_0(g)
\end{eqnarray}
    }
because $(-1)^{\theta(g)}Ad_{\exp(\pi e_3)}(\rho_0(g))=\rho_0(g)$ and
$Ad_{\exp(\pi e_3)}$ changes the sign of $e_1$ and $e_2$ but
preserves $e_3$. Comparing equality~(\ref{eqn:nu}) with
\[
 \nu(\varrho_{w})(g)=\varrho_{-w}(g)= \exp(f_{g}(-w)+h_g(-w))\rho_0(g)
\]
it follows that $f_g$ is odd and $h_g$ even, as claimed. \qed

 \proof[Proof of Lemma~\ref{lem:S}]
 We choose an element $g_0\in\ker(\theta)$ such that
    $
    \rho_0(g_0)=\exp(\alpha_0 e_3),
    $
for some $\alpha_0\in\mathbb R-2\pi\mathbb Z$. We define:
\begin{displaymath}
    \mathcal S=\left\{ \rho\in R(M,SL_2(\mathbb C))
    \,\left\vert\begin{array}{c}
     \rho(m)=\exp(\alpha e_1),
    \rho(g_0)=\exp(\beta_1 e_1+\beta_3 e_3),\\ \textrm{ with }
    \alpha,\beta_1,\beta_3\in\mathbb C
    \end{array}\right.\right\}
\end{displaymath}
The projection $t\co R(M,SL_2(\mathbb C))\to X(M,SL_2(\mathbb C))$
restricts to a map $t\vert_{\mathcal S}\co \mathcal S\to
X(M,SL_2(\mathbb C))$.

Let $e_0$ denote the identity matrix of size $2\times 2$, so that
$\{e_0,e_1,e_2,e_3\}$ is a basis for $M_2(\mathbb C)$ as $\mathbb
C$--vector space. For every $\rho\in R(M,SL_2(\mathbb C))$ and every
$g\in\pi_1M$ we write:
\[
\rho(g)=\sum_{i=0}^3\pi_{i,g}(\rho) e_i.
\]
If we define $F\co R(M,SL_2(\mathbb C))\to\mathbb C^3$ as
    $
     F=(\pi_{2,m},\pi_{3,m},\pi_{2,g_0})
    $,
    then $\mathcal S=F^{-1}(0)$.
    An easy computation shows that the differential of $F$ at
    $\rho_0$ maps $B^1(M,sl_2(\mathbb C))$ isomorphically onto $\mathbb
    C^3$. It follows that $\rho_0$ is a smooth point of $\mathcal S$ and
     that $t\vert_{\mathcal S}$ is locally
    bianalytic. This proves assertions (i) and (ii) of the
    proposition.

If in the construction of $\mathcal S$ we replace $SL_2(\mathbb C)$
by $SU(2)$, then we obtain $\mathcal S'$ and the same construction as
above applies to prove assertion (iii) of the proposition. Finally
assertions (iv) and (v) follow from construction. \qed

\subsection{Sections for the deformation spaces}

\begin{defn} For $(s,t)\in U$, we define $\rho_{(s,t)}\in R(M,G))$ as follows:
\[
\rho_{(s,t)}=\left\{
    \begin{array}{cl}
      \sigma(s-t\, i)\in R(M,SL_2(\C)) & \textrm{ when }\mathbb X^3 =\mathbb H^3 \\
       (\sigma(s+t),\sigma(s-t))\in R(M,SU(2)\times SU(2))
        & \textrm{ when }\mathbb X^3 =\mathbb S^3
    \end{array}
\right.
\]
\end{defn}

\begin{prop}\label{prop:propertyx0} For every $(s,0)\in U$,
$\rho_{(s,0)}$ stabilizes $x_0\in \mathbb X^3$, the  point
stabilized by $\rho_0$.
\end{prop}

\proof Let $f_g$ and $h_g$ be the functions of
Lemma~\ref{lem:sigma0}. In the hyperbolic case, the proposition
follows from the fact that the functions $f_g$ and $h_g$ have real
coefficients: when $t=0$, $f_g(s), h_g(s)\in su(2)$, hence
$\rho_{(s,0)}\in R(M,SU(2))$, and $SU(2)$ is precisely the
stabilizer of $x_0$ . In the spherical case, $\rho_{(s,0)}$ is
diagonal by construction, and the diagonal is precisely the
stabilizer of $x_0$. \qed

\section{Infinitesimal deformations}\label{sect:inf}

\subsection{Infinitesimal isometries}

Recall that in the convention after Theorem~\ref{thm:defstr}, we have
fixed a point $x_0$ so that $\phi_0=\ROT\circ\pi\circ\hol$ is a
representation into
\[
SO(3)\cong \Isom^+(X^3)_{x_0}\hookrightarrow\Isom^+(X^3).
\]
Its lift to $SU(2)\cong G_{x_0}$ is $\rho_0$.
 Let $\mathbf
g$ denote the lie algebra of $\Isom^+(\mathbb X^3)$ and $\mathbf
g_{x_0}$ the Lie subalgebra corresponding to $ \Isom^+(\mathbb
X^3)_{x_0}$. We have a natural exact sequence
\[
0\to \mathbf g_{x_0}\to \mathbf g\to T_{x_0}\mathbb X^3\to 0
\]
The Killing form  on $\mathbf g$ is non-degenerate, and
$T_{x_0}\mathbb X^3$ is naturally identified to the orthogonal space
to $\mathbf g_{x_0}$. We have an orthogonal sum:
\begin{equation}\label{eqn:rotplustrans}
    \mathbf g=\mathbf g_{x_0}\bot T_{x_0}\mathbb X^3
\end{equation}

\begin{defn}
 Elements of $\mathbf g$ are called infinitesimal isometries;
elements of $\mathbf g_{x_0}$, infinitesimal rotations (with respect
to $x_0$); and elements of  $T_{x_0}\mathbb X^3$, infinitesimal
translations (with respect to $x_0$).
\end{defn}

\begin{lem} \label{lem:identification} There is a natural identification of $SO(3)$--modules:
\[
   T_{x_0}\mathbb X^3\cong\R^3\cong \mathbf g_{x_0},
\]
where the action of $SO(3)$ on $\mathbf g_{x_0}$ and
$T_{x_0}\mathbb X^3$ is the adjoint action and the action on
$\R^3$ is standard. In addition, it preserves the products (cross
product on $\mathbb R^3\cong T_{x_0}\mathbb X^3$ and Lie bracket
on $\mathbf g_{x_0}$) and the natural bilinear forms (Killing form
on $\mathbf g_{x_0}$ and the metric on $\mathbb R^3\cong
T_{x_0}\mathbb X^3$) up to a constant. \qed
\end{lem}

 The isomorphism from $ T_{x_0}\mathbb X^3$ to $\mathbf
g_{x_0}$ maps the infinitesimal translation of tangent vector $v\in
 T_{x_0}\mathbb X^3$ to the infinitesimal rotation around the line
 $\R v$ of infinitesimal angle $\vert v\vert$.


It is convenient to specify Lemma~\ref{lem:identification} and
isomorphism~(\ref{eqn:rotplustrans}) in the hyperbolic and the
spherical case:
\begin{enumerate}[(a)]
\item In the hyperbolic case $\mathbf g\cong sl_2(\C)$ and $\mathbf
g_{x_0}$ is a subalgebra conjugate to $su(2)$. In this case,
isomorphism~(\ref{eqn:rotplustrans}) is written as:
\[
sl_2(\C)= \mathbf g_{x_0}\bot i \mathbf g_{x_0}.
\]
In addition the isomorphism of Lemma~\ref{lem:identification} maps
$v\in\mathbf g_{x_0}$ to $-i v\in T_{x_0}\mathbb X^3$.
\item
In the spherical case $\mathbf g\cong su(2)\times su(2)$. Up to
conjugation, $\mathbf g_{x_0}\cong su(2)$ is the subalgebra  of
diagonal matrices and $T_{x_0}\mathbb X^3$ is the set of
anti-diagonal elements (ie, matrices of the form $(a,-a)$ with $a\in
su(2)$). Hence isomorphism~(\ref{eqn:rotplustrans}) is the
decomposition of matrices of $su(2)\times su(2)$ as the sum of
diagonal plus anti-diagonal elements. The isomorphism of
Lemma~\ref{lem:identification} maps $(a,a)\in\mathbf g_{x_0}$ to
$(a,-a)\in T_{x_0}\mathbb X^3$.
\end{enumerate}

As an application we obtain:

\begin{prop}\label{prop:propertyX2}Let  $\mathbb X^2\subset \mathbb X^3$ denote the geodesic
hyperplane preserved by $\rho_0$ (tangent to $\R^2\times 0$). Then
for every $(0,t)\in U$, $\rho_{(0,t)}$ preserves  $\mathbb X^2$.
\end{prop}

\proof In the hyperbolic case we use the fact that $f_g$ is odd
and $h_g$ is even. Hence $f_g(i\, t)$ is purely imaginary and
$h_g(i\, t)$ is real. Thus $f_g(i\, t)$ is an infinitesimal
translation tangent to $\mathbb X^2$ and $h_g(i\, t)$ is an
infinitesimal rotation around a geodesic perpendicular to $\mathbb
X^2 $. This means that $f_g(i\, t)+h_g(i\, t)$ belongs to the Lie
algebra of the isometry group of $\mathbb X^2$.
 In the
spherical case, $(f_g(s),f_g(-s))=(f_g(s),-f_g(s))$ and
$(h_g(s),h_g(-s))=(h_g(s),h_g(s))$, which also means that these
elements belong to the Lie algebra tangent to the isometry group of
$\mathbb X^2$. \qed

\subsection{Infinitesimal properties of the section}

Let $ \partial_s \rho \co \pi_1M\to\mathbf g$ denote the cocycle
defined by
     \[
   g \mapsto\left.
\partial_s(\rho_{(s,t)}(g)\rho_0(g^{-1}))\right\vert_{(s,t)=0}
\qquad\textrm{ for every } g\in\pi_1M.
\]
We use the equivalent notation for $\partial_t$.

\begin{lem}\label{lem:infprop}\begin{enumerate}[\rm(i)]
\item The cocycle $
    \partial_s \rho$ is valued on infinitesimal rotations.
\item The cocycle $\partial_t \rho$ is valued on infinitesimal translations.
\item   Under the identification of Lemma~\ref{lem:identification},
    $\partial_s \rho=\partial_t \rho$.
In addition, they are valued on the invariant plane
$\R^2\times\{0\}$.
\end{enumerate}
\end{lem}

\proof Assertion (i) follows from
Proposition~\ref{prop:propertyx0}. The remaining assertions follow
easily from construction. For instance, in the hyperbolic case,
$\partial_s \rho=f_g'(0)$ and $\partial_t \rho=-i\ f_g'(0)$,
because $h_g'(0)=0$ ($h_g$ is even).
    In the spherical case, $\partial_s \rho=(f_g'(0),f_g'(0))$ and
    $\partial_t \rho=(f_g'(0),-f_g'(0))$. (See
 the explanation after Lemma~\ref{lem:identification}).\qed

\begin{defn} We define $\partial_{s}\partial_{s}\log\rho$
    to be the chain in $C^1(M,\mathbf g)$ such that
 $\forall g\in\pi_1M$,
    \[
(\partial_{s}\partial_{s}\log\rho)(g)=
\frac{\partial^2\phantom{s}}{\partial
s^2}\log(\rho_{(s,t)}(g)\rho_0(g^{-1}))\vert_{(s,t)=0}.
    \]
We use the same definition  for all other partial derivatives.
\end{defn}

\begin{prop} \label{prop:infeq} There exists a choice of $p\in Nil$
and of the holonomy representation $\hol\co\pi_1\mathcal O\to Nil$
 such that, if $q=\pi(p)$,
 then:
\begin{enumerate}[\rm(i)]
    \item $\partial_s\rho=z_q$,
    \item The cochain $\partial_s\partial_t\log\rho$
     is valued on infinitesimal translations along
      the invariant line $0\times\R$ and equals to $c_p$.
    \item The translational part of $\partial_s\partial_s\log\rho$
     and of $\partial_t\partial_t\log\rho$ vanish.
\end{enumerate}
\end{prop}

\proof  The cocycle $z_q\in Z^1(M,\R^2\times 0)$ represents a
non-zero element in cohomology.
 In
addition, $\partial_s\rho$ is also non-zero in cohomology,
    because $w$ is locally a parametrization.
Since $ H^1(M,\R^2\times 0)\cong\mathbb R$, by composing the holonomy
$\hol$ with an automorphism  of $Nil$ of the form
\[
(x_1,x_2,x_3)\mapsto  (\lambda x_1,\lambda x_2,\lambda^2 x_3),
\qquad\textrm{ for all } (x_1,x_2,x_3)\in Nil,
\]
 we have equality
(i) up to coboundary. The choice of $q$ eliminates the indeterminacy
of the coboundary.

 To prove (ii), since $f''_g(0)=0$, in the hyperbolic case we have
 $(\partial_s\partial_t\log \rho)(g)=- i h_g''(0)$ and in the
 spherical case $(\partial_s\partial_t\log
 \rho)(g)=(h_g''(0),-h_g''(0))$.
 In both cases $(\partial_s\partial_t\log \rho)(g)$ is an
 infinitesimal translation with value $h''_g(0)\in 0\times \R$.

 From the second order terms
in the expression $\varrho_{w}(g_1
g_2)=\varrho_{w}(g_1)\varrho_{w}(g_2)$ we obtain:
\[
h_{g_1}''(0)+Ad_{\rho_0(g_1)}(h_{g_2}''(0))+ [f_{g_1}'(0),
Ad_{\rho_0(g_1)}(f_{g_2}'(0))]
=
h_{g_1g_2}''(0)
\]
(use for instance the Campbell-Hausdorff formula). Hence
\[
\partial_s\rho\cup\partial_s\rho =\delta ( \partial_s\partial_t\log\rho
).
\]
Since $ H^1(M,0\times\R)=0$, we have that $c_p$ equals
$\partial_s\partial_t\log \rho$ up to a coboundary. Again the
indeterminacy of the coboundary is eliminated by choosing
conveniently $p\in\pi^{-1}(q)$.

Finally to prove (iii), in the hyperbolic case
$(\partial_s^2\log\rho)(g)=-(\partial_t^2\log\rho)(g) = h_g''(0) $
and in the spherical case
 $(\partial_s^2\log\rho)(g)=(\partial_t^2\log\rho_0)(g) =
 (h_g''(0),h''_g(0))$. In both cases, these are infinitesimal
 rotations.
 \qed

\subsection{Compatibility with the holonomy}

In this subsection we prove
 property (iii) of Proposition~\ref{prop:defrep}; property (iv) being similar is not
 proved. We want to prove
 that for every $g\in\pi_1M$:
\[
\lim_{\begin{smallmatrix}
  (s,t)\to 0 \\
  s\, t\neq 0
\end{smallmatrix}} \Delta_{(s,t)}^{-1} \circ \rho_{(s,t)}(g)\circ \Delta_{(s,t)} =
\hol(g)
\]
uniformly on compact subsets of $Nil$ for the $\mathcal
C^1$--topology.

\proof We fix $g\in\pi_1M$. We know that
\begin{equation}\label{eqn:exp-1}
\exp^{-1}_{x_0}(\rho_{(s,t)}(g)(\Delta_{(s,t)}(x_1,x_2,x_3)))
\end{equation}
is analytic on $(s,t)$ and on $(x_1,x_2,x_3)$. In addition:
\begin{itemize}
\item[--] the expression~(\ref{eqn:exp-1}) is a multiple of $t$,
 because when $t=0$, $\rho_{(s,0)}(g)$
fixes $x_0$ (Proposition~\ref{prop:propertyx0}), and
\item[--] the coefficient
in $e_3$ of (\ref{eqn:exp-1}) is a multiple of $s\, t$, because when
$s=0$, $\rho_{(0,t)}(g)$ preserves $\mathbb X^2=\exp_{x_0}( \mathbb
R^2\times 0)=\exp_{x_0}\langle e_1,e_2\rangle $
(Proposition~\ref{prop:propertyX2}).
\end{itemize}
 Thus it suffices to compute the first order terms
of~(\ref{eqn:exp-1}). More precisely, we write the
expression~(\ref{eqn:exp-1}) as follows:
\[
f_1(x,(s,t)) e_1+ f_2(x,(s,t)) e_2+f_3(x,(s,t)) e_3
\]
for some analytic functions $f_i$ such that $f_1$ and $f_2$ are
multiples of $t$ and $f_3$ is a multiple of $s\, t$. We want to prove
that
\[
(\partial_t f_1(x,0),\partial_t f_2(x,0),\partial_t\partial_s
f_3(x,0))=\hol(g)(x).
\]
We notice that analyticity implies that the convergence is uniform on
compact subsets for the $\mathcal C^1$--topology.

Corresponding to the basis $\{e_1,e_2,e_3\}$
 for the sub-algebra $su(2)=\mathbf g_x$ (ie, infinitesimal
rotations), there is a basis for the space of infinitesimal
translations $\{w_1,w_2,w_3\}$
 via the isomorphism of Proposition~\ref{prop:infeq}.
We have the following relations up to cyclic permutation of
coefficients:
\[
[e_1,e_2]=e_3,\ [w_1,w_2]=k e_3,\ [e_1,w_2]=w_3,\ [e_1,w_1]=0,
\]
where $k=\pm 1$ is the curvature of $\mathbb X^3$. In addition we
have
 \[
\Delta_{(s,t)}(x_1,x_2,x_3)=\exp(t x_1 w_1+t x_2 w_2+st x_3 w_3)
(x_0).
 \]
If $\hol(g)$ is the multiplication by $(a_1,a_2,a_3)\in Nil$ composed
with $\rho_0(\gamma)$, then by Proposition~\ref{prop:infeq}
\[
\rho_{s,t}(g)=\exp (a_1 ( s e_1+ t w_1) +a_2(s e_2+ t w_2) + a_3 st
w_3 + A) \rho_0(g)
\]
where $A$ are higher
 order terms (of order two multiplying $e_1
,e_2,e_3,w_1,w_2$ and of order three multiplying $w_3$).

 Since
$ \rho_0(g)(\Delta_{(s,t)}(x))=\Delta_{(s,t)}(\rho_0(g)(x)) $, we may
assume that $\rho_0(g)$ is trivial. We use the following notation
\[\begin{array}{rcl}
  R & = & s\,a_1 e_1 +s\,a_2 e_2\\
  T & = & t(a_1 w_1 +a_2 w_2+s\, a_3 w_3)\\
  X & = & t(x_1 w_1 +x_2 w_2+s\, x_3 w_3)
\end{array}
\]
so that $\rho_{(s,t)}(g)=\exp(R+T+A)$ (we are assuming that
$\rho_0(g)=  \operatorname{Id}$) and
$\Delta_{(s,t)}(x_1,x_2,x_3)=\exp(X)(x_0)$. Hence:
    \[
    \rho_{(s,t)}(g)(\Delta_{(s,t)}(x_1,x_2,x_3))=\exp(R+T+A)\exp(X)(x_0).
    \]
By the Campbell-Hausdorff formula:
\begin{eqnarray*}
    \lefteqn{\exp(R+T+A)\exp(X)=     \exp(R+T+X+\frac12[R+T,X]+A)  }\\
    &=& \exp(T+X+\frac12[R+T,X]-\frac12[T+X,R]+A)\exp(R)\\
    &=&\exp(T+X+ [R,X]+A)\exp(R) ,
\end{eqnarray*}
where $A$ is as above, because $[T,X]$ is an infinitesimal
rotation of order two and $[R,T]$ is a translation but of order
three.  Since $\exp(R)(x_0)=x_0$, it follows that
\[
\rho_{(s,t)}(g)(\Delta_{(s,t)}(x_1,x_2,x_3))=\exp(T+X+[R,X]+A)(x_0).
\]
In addition $[R,X]=(a_1 x_2-a_2 x_1)s\, t\, w_3+O(s^2t)$, and
property (iii) of Proposition~\ref{prop:defrep} follows. \qed

 \section{Euclidean structures}\label{sec:euclidean}

In this section we prove the part of Theorem~\ref{thm:defstr}
concerning Euclidean structures. We use the semi-direct product
structure of the isometry group
 and its universal covering:
\[
\Isom^+(\mathbb R^3)\cong \R^3\rtimes SO(3),\qquad
\widetilde{\Isom^+}(\mathbb R^3)\cong\R^3\rtimes SU(2).
\]

\begin{defn}\label{def:Phi}
For $s$ in a neighborhood of the origin, we define the representation
 $\rho'_s\co \pi_1M\to
\R^3\rtimes SU(2)$ as:
\[
\rho'_s=\left( \partial_t\rho_{(s,0)},\rho_{(s,0)}   \right)
\]
\end{defn}

Notice that $\rho'_s$ is a representation because
$\partial_t\rho_{(s,0)}$ is a cocycle twisted by $\rho_{(s,0)}$. In
particular $\ROT\circ\rho'_s=\rho_{(s,0)}$. The action of
$\rho'_s(g)$ on $\mathbb R^3$ is the following:
\[
v\mapsto
\rho_{(s,0)}(g)(v)+\partial_t\rho_{(s,0)}(g)\qquad\qquad\forall
v\in\mathbb R^3\cong T_{x_{0}}\mathbb X^3.
\]

\begin{defn}
 We define the map $D'_s\co \widetilde M\to T_{x_0}\mathbb X^3\cong\R^3$
as
\[
D'_s(x)=\partial_t D_{(s,t)}(x)\vert_{t=0}
\]
\end{defn}

Since $D_{(s,0)}$ is the constant map $x_0$, the image of $D'_s$ is
contained in $T_{x_0}\mathbb X^3$.

The following proposition shows that $D'_s$ is a developing map

\begin{prop}\label{prop:Ds} The map $D'_s$ is $\rho'_s$--equivariant and it is a
local diffeomorphism for $s\neq 0$.
\end{prop}

The proof requires the following lemma.

\begin{lem}\label{lem:dpsi} Let $\gamma\co (-\varepsilon,\varepsilon)\to \mathbb X^3$
be a path such that $\gamma(0)=x_0$ and $\gamma'(0)=v\in
T_{x_0}\mathbb X^3$. Then
\[
\rho'_s(g)(v)=\partial_t(\rho_{(s,t)}(g)(\gamma(t)))\vert_{t=0}.
\]
\end{lem}

\proof By the chain rule:
\begin{multline*}
    \partial_t \rho_{(s,t)}(g)(\gamma(t))\vert_{t=0}=
    \partial_t \rho_{(s,t)}(g)(x_0)\vert_{t=0} +
    \rho_0(g)(\gamma'(0))=\\
    \partial_t \rho_{(s,0)}(g)+\rho_0(g)(v)=\rho'_s(g)(v).
    \qed
\end{multline*}

\proof[Proof of Proposition~\ref{prop:Ds}] Equivariance of $D'_s$
follows from deriving the following equality
\[
\rho_{(s,t)}(g) (D_{(s,t)}(x))= D_{(s,t)}(g\cdot x)
\]
and applying Lemma~\ref{lem:dpsi}.

Next we write $\Delta'_s(x) =\partial_t\Delta_{(s,t)}(x)\vert_{t=0}$.
We have
\[
\Delta'_s(x_1,x_2,x_3)= x_1e_1+x_2e_2+s\,x_3e_3
\]
Hence $\Delta'_s$ is a diffeomorphsim for $s\neq 0$. We claim that
\[
\lim_{s\to 0} (\Delta'_{s})^{-1} \circ \rho'_s(g)\circ
\Delta'_s=\hol(g).
\]
The proof of this claim follows a scheme similar to the proof of
Proposition~\ref{prop:defrep} (iii) and deriving with respect to $t$
some of its equalities.

The construction of $D_{(s,t)}$ implies that $D'_s$ can also be
constructed by using bump functions. Hence
 $ (\Delta_s')^{-1}\circ D'_s$ converges to
 $\mathcal D_0$ uniformly on compact subsets, and
the proposition follows. \qed

\section{Deformation space and Dehn filling coefficients}\label{sect:DF}

In this section we construct the deformation space {\it Def}, we define
the Dehn filling coefficients and we study its behaviour. At the end
of the section we prove Theorem~B.

\begin{defn}
    We define the deformation space {\it Def} as the open set
    \[
        {\it Def}=\{(s,\tau)\in\mathbb R^2\mid s\geq 0, (s,\tau)
        \textrm{ in a neighborhood of } 0\}
    \]
    such that $(s,\tau)$ corresponds to the structure with parameters
    $(s,t)$ as follows:
    \begin{itemize}
        \item[--] when $\tau>0$, it corresponds to the hyperbolic structure with
        $\tau=t^2$,
        \item[--] when $\tau<0$, to the
        spherical structure with
        $\tau=-t^2$,
        \item[--] when $\tau=0$,  to the
        Euclidean  structure with
        $t=0$.
    \end{itemize}
\end{defn}

\subsection{Dehn filling coefficients}

We shall define the Dehn filling coefficients and prove that they
induce  an analytic map on $(s,\tau)\in {\it Def}$:
\[
(p,q)\co {\it Def}\to \mathbb R^2.
\]
We recall that in Subsection~\ref{subsec:localpara} we have chosen
$l,m\in\pi_1M$ that generate a peripheral subgroup $\pi_1T^2$, so
that $m$ is a meridian of $\Sigma$. We notice that since $l$ and $m$
commute, their holonomies have a common invariant geodesic.

\begin{defn}
For a geometric structure with holonomy $\rho_{(s,t)}$, we define
$u\in\mathbb C$ to be the complex length of $\rho_{(s,t)}(m)$ (ie,
$\rho_{(s,t)}(m)$ is translation of length $\re(u)$ composed with a
rotation of angle $\im(u)$ along the invariant geodesic). We also
define $v\in\mathbb C$ as the complex length of $\rho_{(s,t)}(l)$.
 \end{defn}

The parameters $(u,v)$ are not uniquely defined. Besides $(u,v)$ we
could choose any pair in the following set:
\[
\pm (u+ 2\pi i\mathbb Z, v+2\pi i\mathbb Z).
\]
The choice of the sign depends on the orientation of the geodesic
invariant by $\rho_{(s,t)}(l)$  and $\rho_{(s,t)}(m)$.
 We view $u(s,t)$ and $v(s,t)$ as analytic functions
on $(s,t)$, hence they are unique if we fix the branch with
$u(0,0)=\pi i$ and $v(0,0)=0$.

\begin{defn}
    Given $(s,t)\in U$, $(p,q)\in\mathbb R^2$ are defined by the rule
    \[
    p u+q v=2\pi i.
    \]
\end{defn}

This definition is equivalent to:
\begin{equation}\label{eqn:DFcoeff}
    \left.
    \begin{array}{rcl}
        p \re u + q \re v &=&0\\
        p \im u + q \im v &=& 2\pi
    \end{array}
    \right\}
\end{equation}

\begin{prop}
    If we fix the branch $p(0,0)=2$ and $q(0,0)=0$, then
    $(p,q)$ is an analytic map on $(s,\tau)\in Def$.
\end{prop}

\proof We start by describing $(u,v)$ as analytic maps on $(s,t)$ in
the hyperbolic, spherical and Euclidean cases.

Let $w$ be the local parameter of Definition~\ref{def:w}, and let
$\varrho_{w}=\sigma(w)$, where $\sigma$ is the section in
Lemma~\ref{lem:sigma0}. By Lemmas~\ref{lem:sigma0} and~\ref{lem:S}
(v), there exists an odd analytic function $F$ with real coefficients
such that
\[
\varrho_{w}(m)=\pm\exp((\pi+F(w)) e_1)=\pm \exp(F(w) e_1)\exp(\pi
e_1).
\]
Since $m$ and $l$ commute, by definition of $w$ we have:
\[
\varrho_{w}(l)=\pm\exp(w\, e_1).
\]
In the hyperbolic case, $w=s-i\, t$ (see
Subsection~\ref{subsec:defchar}), hence:
\[
    \left\{
        \begin{array}{l}
            u_H= i(\pi+F(s-i\, t))=\im (F(s+i\, t))+ i \, (\pi+\re(F(s+i\, t)))\\
            v_H= i(s-i\, t)= t+i\, s
        \end{array}
    \right.
\]
In the spherical case, we work in $X(M,SU(2))\times X(M,SU(2))$ and
we take $(w_1,w_2)=(s+t,s-t)$ (see also
Subsection~\ref{subsec:defchar}). Hence:
\[
    \left\{
        \begin{array}{l}
            u_S= (F(s+ t)-F(s-t))/2+ i \, (\pi+(F(s+ t)+F(s-t))/2)\\
            v_S= t+i\, s
        \end{array}
    \right.
\]
In the Euclidean case the translational part is obtained by deriving
with respect to $t$ when $t=0$ (see Section~\ref{sec:euclidean}).
Thus:
\[
    \left\{
        \begin{array}{l}
            u_E= F'(s)+ i \, (\pi+F(s))\\
            v_E= 1+i\, s
        \end{array}
    \right.
\]
Before showing that $(p,q)$ are well defined, we must notice that
$\re(u_H)$ and $\re(u_S)$ are both multiples of
$t=\re(v_H)=\re(v_S)$. Hence we redefine:
\[
    \left\{
        \begin{array}{l}
            \tilde u_H= \im (F(s+i\, t))/t+ i \, (\pi+\re(F(s+i\,
            t)))\\
            \tilde v_H= 1+i\, s\\
           \tilde  u_S= (F(s+ t)-F(s-t))/(2t)+ i \, (\pi+(F(s+
           t)+F(s-t))/2)\\
            \tilde v_S= 1+i\, s.
        \end{array}
    \right.
\]
We keep $\tilde u_E=u_E$ and $\tilde v_E=v_E$. The system of
equations~(\ref{eqn:DFcoeff}) becomes
\begin{equation}
    \left.
    \begin{array}{rcl}
        p \re \tilde u + q \re \tilde v &=&0\\
        p \im \tilde u + q \im \tilde v &=& 2\pi
    \end{array}
    \right\}
\end{equation}
Since $\re \tilde v=1$, $\im\tilde u=\pi+O(s,t)$, $\re\tilde
u=O(s,t)$ and $\im\tilde v=s$, it is clear from this system of
equations that $(p,q)$ is a well-defined analytic map on $(s,t)$ in
every case (hyperbolic, Euclidean and spherical).

To show that $(p,q)$ is an analytic map on $(s,\tau)\in {\it Def}$, we must
check the following properties:
\begin{enumerate}[(i)]
    \item $\tilde u_H(s,t)=\tilde u_S(s,i\, t)$.
    \item $\tilde u_H(s,0)=\tilde u_S(s,0)=\tilde u_E(s)$
    \item $\tilde u_H(s,t)$ and $\tilde u_S(s, t)$ are even on $t$.
\end{enumerate}
These  properties are obvious from construction. \qed

\subsection{The power expansion of $(p,q)$}

In this section we compute the power expansion of $(p,q)$.
First we
need the following proposition.

\begin{prop}\label{prop:Fa3} $F(w)=a_3 w^3+O(w^5)$, with $a_3>0$.
\end{prop}

\begin{lem}\label{lem:F'=0} $F'(0)=0$. \end{lem}

\proof Using the notation of Lemma~\ref{lem:w}, $\alpha_m=\pi+F(w)$.
In the same lemma it is proved that $d\alpha_m=0$, thus
$F'(0)=0$.\qed

\proof[Proof of Proposition~\ref{prop:Fa3}]
 We know that $F$ is an odd function with $F'(0)=0$. In the
proof we use Theorem~\ref{thm:defstr}: there is a neigbhorhood
$U\subset \mathbb R^2$ of the origin  such that for every $(s,t)\in
U$ with $s\, t\neq 0$, $\rho_{(s,t)}$ is the holonomy of a hyperbolic
structure on $M$ with end of Dehn filling type. The structure at the
end is described by $u$ and $v$.

 We first
show that $F$ is not constant by contradiction.
 If $F$ is constant,
then $F\equiv 0$  because $F$ is odd, and $u\equiv \pi i$. This
implies that all the structures on $U$ induce hyperbolic cone
structures with cone angle $\pi$. This is impossible, because it
implies that $\OO$ is hyperbolic.

Let $2 n+ 1\geq 3$ be the order of the first derivative such that
$F^{(2n+1)}(0)\neq 0$. We claim that $2 n+1=3$. Identifying $\mathbb
C\cong\mathbb R^2$ via $w=s-t\, i$, the map $F\vert_{U}$ is a
branched covering of the open set $F(U)\subset\mathbb C$. It is
branched at the origin with branching order $2 n+1$. We look at the
inverse image of the real line $(F\vert_{U})^{-1}(\mathbb R)$. It
consists of $2n+1$ curves passing through the origin. One of them is
real, hence it corresponds to $t=0$ in $U$. The other $2 n$ curves,
are contained in $\{(s,t)\in U\mid s\, t\neq 0 \}$, hence they give
geometric structures, except for the origin. Since the image of these
curves is real, they correspond to cone structures.

 The intersection of these $2 n$
curves with $\{(s,t)\in U\mid s\, t\neq 0 \}$ has  $4 n$ components,
(each curve is divided into two when we remove the origin). Thus
there are $n$ curves on the quadrant $\{(s,t)\in U\mid s>0, t> 0 \}$.
If $n\geq 2$, then there would be at least two curves in the same
quadrant. These two curves correspond to two families of structures
with the same orientation and spin structure. In addition, when we
parameter the curves from the origin, one of them has decreasing cone
angle $\alpha_m=\pi+F(w)$, but the other one has increasing cone
angle. This is not possible, because Schl\"afli's formula implies
that the cone angles must decrease. Hence $n=1$ and $2 n+1=3$.

Finally, the argument above gives $a_3>0$, because the branch of the
first quadrant corresponds to decreasing volume. \endproof

We will determine the power expansion of $(p,q)$ by analyzing its
behavior on the curves $s=0$, $\tau=-s^2$ and $\tau=0$.

\begin{lem}\label{lem:s=0}
    The Dehn filling coefficients $(p,q)$ induce a bijection between the
    segments $s=0$ and $p=2$;
\end{lem}

\proof  Since $F$ is odd, $F(i\, t)$ has zero real part. Hence,
when $s=0$, $\tilde u_H(0,t)=F(i\ t)/(i\, t)+i\,\pi$ and $\tilde
v_H(0,t)=1$. Therefore in the hyperbolic case
\begin{equation}\label{eqn:s=0hyp}
    p=2,\quad q=-2 F(i\,t)/(i\, t)=2 a_3 t^2+O(t^4)=2a_3 \tau
    +O(\tau^2).
\end{equation}
In addition, $\tilde u_S(0,t)=F(t)/t+i\pi$ and $\tilde
v_S(0,t)=1$. Thus in the spherical case
\begin{equation*}
    p=2,\quad q=-2F(t)/t=-2a_3 t^2+O(t^4)=2a_3\tau+O(\tau^2),
\end{equation*}
which is the same as equation~(\ref{eqn:s=0hyp})   but for
$\tau<0$. \qed

\begin{lem}\label{lem:tau=-s^2}
The Dehn filling coefficients $(p,q)$ induce a bijection between the
curve $\tau=-s^2$ and the segment $p=2$, $q\leq 0$.
\end{lem}

\proof Since the equation $\tau=-s^2$ is equivalent to $t=s$ in
the spherical case, we have
    $ \tilde u_S(s,s)=F(2s)/(2s)+ i (\pi+
    F(2s)/2)
    $
and $v_S(s,s)=1+i\, s$. This gives the curve:
\begin{equation}\label{eqn:tau=-s^2}
        p=2,\quad
        q=-F(2s)/s= -8 a_3 s^2+O(s^4) =8 a_3\tau+O(\tau^2)
\end{equation}
Hence the lemma is clear. \qed

\begin{rem}
 The structures of Lemma~\ref{lem:s=0} are transversely
riemannian foliations. The structures of Lemma~\ref{lem:tau=-s^2}
are spherical and they are equipped with an isometric foliation of
codimension 2 (in particular it is also transversely spherical).
This comes from the fact that the equation $s=t$ implies that the
parameter in Subsection~\ref{subsec:defchar} is
$(w_1,w_2)=(2s,0)$. Hence the image
 of the holonomy representation is contained in $SU(2)\times \widetilde{O(2)}$,
where $\widetilde{O(2)}$ is the lift of $O(2)< SO(3)$. Hence it is
compatible with the isometric action of $\{1\}\times S^1<SU(2)\times
SU(2)$.
\end{rem}

\begin{lem}\label{lem:t=0}
    The Dehn filling coefficients map the half line $\tau=0$ bijectively to a half
    curve with power expansion:
    \[
        \left\{\begin{array}{l}
        p=2+\frac{4a_3}{\pi}s^3+O(s^5)\\
        q=-6 a_3 s^2+O(s^4)
        \end{array}
        \right.
    \]
\end{lem}

\proof When $\tau=0$, $u_E=F'(s)+i(\pi+F(s))$ and $v_E=1+i\, s$.
Hence
\[
p=2\pi/(\pi+F(s)-s\, F'(s))\quad\textrm{ and } \quad q=-p \, F'(s).
\]
Since $F(s)=a_3 s^3+ O(s^5)$, the lemma is straightforward. \qed

\begin{defn} We define
$g:(-\varepsilon,\varepsilon)\to\mathbb R$ to be a real function such
that, for $q\geq 0$, $g(q)=2$, and for $q\leq 0$, $p=g(q)$ is the
half curve of Lemma~\ref{lem:t=0}.
\end{defn}

\begin{cor}\label{cor:powerexpansion}
    We have the following power expansion:
    \[
        \left\{\begin{array}{l}
        p=2+s(s^2+\tau)(\frac{4a_3}{\pi}+O(s,\tau))\\
        q=2 a_3 (\tau-3 s^2)+O(\tau s^2)+O(\tau^2)+O(s^3)
        \end{array}
        \right.
    \]
\end{cor}

\proof By Lemmas~\ref{lem:s=0} and \ref{lem:tau=-s^2}, $p-2$ is a
multiple of $s(\tau+s^2)$. The coefficient $\frac{4a_3}{\pi}$ comes
from Lemma~\ref{lem:t=0}. The power expansion of $q$ is
straightforward from equation~(\ref{eqn:s=0hyp}) and
Lemma~\ref{lem:t=0}. We notice that $q$ has no coefficient in $\tau
s$, by equation~(\ref{eqn:tau=-s^2}). \qed

\subsection{The Whitney pleat}\label{subsec:Whitney}

In the next proposition we view $(p,q)$ as a function on $(s,\tau)$
defined not only on {\it Def} but in a neighborhood of the origin in
$\mathbb R^2$.

\begin{prop}\label{prop:Whithney}
    The map $(p(s,\tau),q(s,\tau))$ has a Whitney pleat at the
    origin, with folding curve $\tau=-9 s^2+O(s^3)$.
\end{prop}

\proof Using the power expansion of
Corollary~\ref{cor:powerexpansion}, the Jacobian is:
\[
J(s,\tau)=\left\vert
    \begin{array}{cc}
    p_s&p_{\tau}\\
    q_s&q_{\tau}
    \end{array}
    \right\vert =
    \frac{8 a_3^2}{\pi}(9 s^2+\tau)+O(\tau s^2)+O(\tau^2)+O(s^3)
\]
Hence $J=0$ is a curve with power expansion $\tau=-9 s^2+O(s^3)$. To
show that there is a Whitney pleat with folding curve $J=0$, we
compute the power expansion of $q$ restricted to this curve:
\[
\phi(s)=q(s,-9 s^2+O(s^3))=-24 a_3 s^2+O(s^3).
\]
Since $\phi''(0)=-48 a_3\neq 0$, the proposition follows
\cite{Whitney}. \endproof

The image of the folding curve $J=0$ is a curve with power expansion:
\begin{equation}\label{eqn:f}
\left\{
    \begin{array}{l}
    p=2-\frac{32 a_3}{\pi} s^3+O(s^4)\\
    q=-\frac{24 a_3}{\pi} s^2+ O(s^3)
    \end{array}
\right.
\end{equation}

\begin{defn} We define
$f\co (-\varepsilon,\varepsilon)\to\mathbb R$ to be a real function
such that, for $q\geq 0$, $f(q)=2$, and for $q\leq 0$, $p=f(q)$ is
the image of the folding curve $J=0$, with $s\geq 0$.
\end{defn}

\begin{figure}[htbp]
 \def\epsfsize#1#2{1#1}
 \centerline{\epsfbox{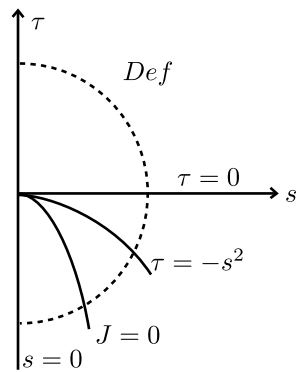}} \caption{
 The curves in the proof of Theorem~B.
 The folding curve is $J=0$, and it is mapped
to $p=f(q)$. The curves $s=0$ and $\tau=-s^2$ are mapped to $p=2$.
The segment $\tau=0$ is the Euclidean region, and it is mapped to
$p=g(q)$}
 \protect \label{fig:Curves}
 \end{figure}

\proof[Proof of Theorem~B] It is clear from
Proposition~\ref{prop:Whithney} and Lemmas~\ref{lem:s=0},
\ref{lem:tau=-s^2} and \ref{lem:t=0}. Notice that the restriction
of $(p,q)$ to {\it Def} gives only half of the Whitney pleat, as in
Figure~\ref{fig:Whitney}. The curves that relevant in the proof of
Theorem~B are recalled in Figure~\ref{fig:Curves}.\qed

\section{The path of cone structures}\label{sect:path}

In this section we prove Propositions~\ref{prop:properties1} and
\ref{prop:properties2} by using the path of cone manifolds. We also
prove the last statement of Theorem~A concerning the limit when
rescaling  those cone manifolds.

Cone structures are determined by the equality $q=0$. From the
power expansion of Corollary~\ref{cor:powerexpansion}, it is clear
that $q=0$ defines a curve in {\it Def}. This curve can be
parametrized as:
\[
\tau= 3 s^2+ O(s^3).
\]
Since $\tau>0$ those structures are hyperbolic. The other coefficient
is $p=2+\frac{16}{\pi}a_3s^3+O(s^4)$. Thus the cone angle is:
\[
\alpha={2\pi}/p=\pi-8a_3s^3+O(s^4)
\]
and therefore the path of cone structures is:
\[\left\{
\begin{array}{rcl}
    s(\alpha) &=&
    \frac12 \sqrt[3]{\frac{\pi-\alpha}{a_3}}
    +O(\vert\pi-\alpha\vert^{2/3})\\
    t(\alpha)&=&\frac{\sqrt{3}}2\sqrt[3]{\frac{\pi-\alpha}{a_3}}
    +O(\vert\pi-\alpha\vert^{2/3}).
\end{array}
\right.
\]
Next we compute some magnitudes of those cone manifolds using the
parameter $s$. The length of the singular set is
\[
\operatorname{length}(\Sigma_{\alpha})=\re(v)=t=\sqrt{\tau}=\sqrt{3}s+O(s^2).
\]
Thus, by Schl\"afli's formula the variation of volume is
\[
d\operatorname{vol}(C_{\alpha})=-\frac12
\operatorname{length}(\Sigma_{\alpha}) d\alpha= (12\sqrt3
a_3s^3+O(4))ds.
\]
Therefore
\[
\operatorname{vol}(C_{\alpha})=3\sqrt 3 a_3s^4+O(5).
\]

\proof[Proof of Proposition~\ref{prop:properties1}]
Straightforward  from the computations above.\hfill$\sq$\medskip

 Below we use that
\[
l_0=\lim_{\pi\to\alpha}\frac{\operatorname{length}(\Sigma_{\alpha})}{(\pi
-\alpha )^{1/3}}= \frac{\sqrt3}{2 a_3^{1/3}} . \]

\proof[Proof of Proposition~\ref{prop:properties2}] We use
the descriptions of the curves $p=f(q)$ and $p=g(q)$ when $q<0$ given
in previous section. First at all, the parametrization of $p=f(q)$
when $q<0$ has a power expansion described in equation~(\ref{eqn:f})
(Subsection~\ref{subsec:Whitney}). Therefore:
\[
\lim_{q\to 0^-}\frac{2-f(q)}{\vert q\vert ^{3/2}}= \frac{\sqrt 2}
{3\sqrt3\pi\sqrt{a_3}}= \frac{4}{9\sqrt[4] 3\pi}l_0^{3/2}
\]
The curve  $p=g(q)$ has a power expansion described in
Lemma~\ref{lem:t=0}, when $q<0$. Thus:
\[
\lim_{q\to 0^-}\frac{g(q)-2}{\vert q\vert ^{3/2}}=
\frac{\sqrt2}{3\sqrt3\pi\sqrt{a_3}}= \frac{4}{9\sqrt[4]
3\pi}l_0^{3/2},
\]
which proves Proposition~\ref{prop:properties2}.\endproof

The following proposition finishes the proof of Theorem~A.

\begin{prop}
When $\alpha\to\pi^-$, the cone manifolds $C_{\alpha}$ re-scaled
by $(\pi-\alpha)^{-1/3}$ converge to the orbifold basis of the
Seifert fibration of $\mathcal O$. In addition, when they are
re-scaled by $(\pi-\alpha)^{-1/3}$ in the horizontal direction and
by $(\pi-\alpha)^{-2/3}$ in the vertical one, they converge to
$\OO$.
\end{prop}

\proof Let $\pi\co Nil\to\mathbb R^2$ denote the projection of the
Riemannian fibration of $Nil$, ie, $\pi(x_1,x_2,x_3)=(x_1,x_2)$. The
developing map of the transverse structure of the Seifert fibration
of $\OO$ is
\[
\pi\circ\mathcal D_0\co\widetilde\OO\to\mathbb R^2
\]
where $\mathcal D_0\co\widetilde\OO\to Nil$ is the developing map
of the $Nil$--structure.

Let $(s(\alpha),t(\alpha))$ denote the path of cone structures.
Since $t(\alpha)$ has order $(\pi-\alpha)^{-1/3}$, to prove the
first part of the proposition is sufficient to show that
\[
\lim_{\alpha\to\pi^-}\frac1{t(\alpha)}\exp_{x_0}^{-1}\circ
D_{(s(\alpha),t(\alpha))} =\pi\circ D_0
\]
uniformly on compact subsets of $\widetilde M$. To prove this limit,
we write
\[
\frac1{t(\alpha)}\exp_{x_0}^{-1}\circ D_{(s(\alpha),t(\alpha))}
=\frac1{t(\alpha)}\exp_{x_0}^{-1}\circ\Delta_{(s(\alpha),t(\alpha))}
\circ \Delta_{(s(\alpha),t(\alpha))}^{-1} \circ
D_{(s(\alpha),t(\alpha))}.
\]
By the proof of Theorem~\ref{thm:defstr},
$\Delta_{(s(\alpha),t(\alpha))}^{-1} \circ
D_{(s(\alpha),t(\alpha))}\to \mathcal D_0$. In addition
\[
\frac1{t(\alpha)}\exp_{x_0}^{-1}\circ\Delta_{(s(\alpha),t(\alpha))}
(x_1,x_2,x_3)=(x_1,x_2,s(\alpha)x_3).
\]
Since $s(\alpha)\to 0$, the limit is clear. Notice that since
$s(\alpha)$ has also order $(\pi-\alpha)^{-1/3}$, the second part
of the proposition follows easily. \qed

\section{An example}\label{sec:example}

\begin{figure}[htbp]
 \def\epsfsize#1#2{1#1}
 \centerline{\epsfbox{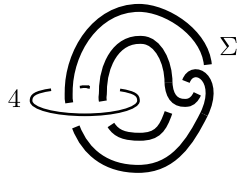}} \caption{
 The surgery description of the orbifold $\OO$.}
 \protect \label{fig:example}
 \end{figure}

We consider the orbifold $\OO$ described as follows. Its
underlying space is the lens space $L(4,1)$, which we view as the
result of Dehn surgery on the trivial knot in $S^3$ with surgery
coefficient 4. We view this trivial knot as one component of the
Whitehead link, and the branching locus $\Sigma$ is precisely the
other component of the link (see Figure~\ref{fig:example}).

It is well known that the Whitehead link has a Montesinos
fibration. This induces an orbifold Seifert fibration of $\OO$. By
looking at the basis of this fibration and its Euler number, one
can check that $\OO$ has $Nil$ geometry.

We want to compute the limit $l_0$ of
Propositions~\ref{prop:properties1} and \ref{prop:properties2}. To
do it, we consider the variety of characters of $M=\OO-\Sigma$.
The manifold $M$ is a punctured torus bundle over the circle, with
homological monodromy
 $
 \left(\begin{smallmatrix}
 1 & 4\\ 1& 5
 \end{smallmatrix}\right)
 $.
Thus its fundamental group admits a presentation
\[
\pi_1(\OO-\Sigma)=\langle a,b,{m}\mid{m} a{\, m}^{-1}=ab,
    {m}\, b\, {m}^{-1}=b (ab)^4
 \rangle
\]
where ${m}$ is the meridian of the branching locus. We also choose
${l}=aba^{-1}b^{-1}$, so that ${l},{m}$ generate a peripheral
group. The variety of characters can be easily computed by using
the methods of~\cite{PoOne}. To compute $l_0$, we do not need the
whole variety of characters, but only its projection to the plane
generated by the variables $x=I_{{m}}$ and $y=I_{{l}}$. This
projection can be computed by means of resultants and it gives the
planar curve:
\[
(y-2)^3+x^2\big(64-16 x^2+x^4+(y-2)(32-5
x^2)+(y-2)^2(7-5y^2)\big)=0
\]
The projection of $\chi_0$ to this curve has coordinates
$(x,y)=(0,2)$.

Using the results of Section~\ref{sect:DF}, we write
\[
y=2\cosh(iw/2)\qquad\textrm{ and }\qquad
x=2\cosh\big(i(\pi+F(w))/2\big).
\]
Since $F(w)=a_3 w^3+O(w^5)$, we have that
\[
y=2-w^2/2+O(w^4)\qquad\textrm{ and }\qquad x=-a_3 w^3+O(w^5).
\]
By replacing those values in the
 the equation of the curve above we obtain:
\[
-(w/2)^6+(a_3 w^3)^2 64 +O(w^8)=0.
\]
Hence $a_3=2^{-6}$. Since $l_0= {\sqrt3}/({2 a_3^{1/3}})$, this
implies that
\[
\lim_{\alpha\to\pi^-}
\frac{\operatorname{length}(\Sigma_{\alpha})}{
(\pi-\alpha)^{1/3}}=l_0=2\sqrt 3 \quad \textrm{ and } \quad \lim
\limits_{q\to
    0^-}  \frac{2-f(q)}{\vert q \vert
    ^{3/2}}=\frac{8\sqrt{2}}{3\sqrt3\pi}.
\]

\section{Cohomology computations}\label{sect:cohmology}

The aim of this section is to prove:
\[
 H^1(M,\R^2\times 0)\cong \R\qquad\textrm{ and } \qquad
H^1(M,0\times \R)=0 .
\]
First we need to compute the homology of the orbifold $\OO$, that can
be defined as follows. Let $K$ be a triangulation of the underlying
space of $\OO$ compatible with $\Sigma$. It induces a triangulation
$\tilde K$ of $\tilde \OO\cong Nil$.
 Let $V$ be a $\pi_1\OO$--module. We consider the following chain
and cochain complexes:
\[
\begin{array}{rl}
 C_*(K;V)&= V \otimes_{\pi_1\OO} C_*(\tilde K;\mathbb Z)
\\ C^*(K;V)&= \Hom_{\pi_1\OO}(C_*(\tilde K;\mathbb Z),V)
\end{array}
\] The homology of $C_*(K;V)$ is denoted by $H_*(\OO;V)$ and the
cohomology of $C^*(K;V)$ by $H^*(\OO;V)$. From the differential point
of view, $H^*(\OO;V)$ is the cohomology of the $V$--valued
differential forms on $\widetilde{\OO}\cong Nil$ which are
$\pi_1\OO$--equi\-variant. The same construction holds for $\Sigma$
and for a tubular neighborhood $\NN(\Sigma)$.

We shall apply the Mayer Vietoris exact sequence to the pair
$(M,\NN(\Sigma))$, so that $M\cup\NN(\Sigma)=\OO$. We first compute
the cohomology of $\OO$.

\begin{lem} Let $V$ be either $\mathbb R^2\times 0$ or $0\times\mathbb
R$. There is a natural isomorphism $H^*(\OO,V)\cong H^*(\pi_1\OO,V)$.
\end{lem}

\proof Let $P\to\OO$ be a finite regular covering such that $P$ is
a manifold. Let ${\Gamma}$ be the group of deck transformations of
the covering.
There is a natural isomorphism
\[
H^*(\OO,V)\cong H^*(P,V)^{\Gamma}.
\]
(See \cite{Bre} for instance). We also have
 a natural
isomorphism
\[
H^*(\pi_1\OO,V)\cong H^*(\pi_1P,V)^{\Gamma}.
\]
 Since $P$ is an
aspherical manifold, there is another natural isomorphism
\[
H^*(\pi_1P,V)^{\Gamma}\cong H^*(P,V)^{\Gamma}.
\]
Hence the lemma follows by composing the three isomorphisms. Notice
that since $C_*(\tilde K;\mathbb Z)$ is an acyclic $\pi_1\OO$--module,
there is a natural map $H^*(\pi_1\OO,V)\to H^*(\OO,V)$, by homology
theory, and that it is the composition of the three isomorphisms.
\qed

\begin{lem}\label{lem:HRx0} $H^0(\OO,\mathbb R^2\times 0)\cong 0$ and
 $H^1(\OO,\mathbb R^2\times 0)\cong \mathbb R$.
 \end{lem}

\proof Since $H^0(\OO,\mathbb R^2\times 0)\cong H^0(\pi_1\OO, \mathbb
R^2\times 0)\cong (\mathbb R^2\times 0)^{\pi_1\OO}$, this group is
zero because the unique element of $\mathbb R^2\times 0$ invariant by
$\pi_1\OO$ is zero.

To compute $H^1(\OO,\mathbb R^2\times 0)$ we use the regular
covering $P\to \OO$ of the previous proof, with deck
transformation group ${\Gamma}$, and the isomorphism
$H^1(P,\mathbb R^2\times 0)^{\Gamma}\cong H^1(\OO,\mathbb
R^2\times 0) $.
 Since the image of $\phi_0$ is finite, we may assume that
 $\pi_1P<\ker\phi_0$. Hence the action of $\pi_1P$ on $\mathbb R^2\times 0$ is
trivial and
\[
H^*(P,\mathbb R^2\times 0)\cong \Hom(H_*(P,\mathbb R), \mathbb
R^2\times 0).
\]
The manifold $P$ can be assumed to be a $S^1$--bundle over $T^2$
with non-trivial Euler number $e\neq 0$. In particular,
\[
\pi_1P\cong\langle t,\alpha,\beta\mid [t,\alpha]=[t,\beta]=1,
[\alpha,\beta]=t^e\rangle
\]
Thus the projection $P\to T^2$ induces a  isomorphism
 $H_1(P,\mathbb R)\cong H_1(T^2,\mathbb R)$
 and:
\[
H^1(P,\mathbb R^2\times 0)\cong\Hom(H_1(T^2,\mathbb R),\mathbb
R^2\times 0)\cong M_{2\times 2}(\mathbb R).
\]where $ M_{2\times 2}(\mathbb R)$
denotes the ring of $2\times 2$ matrices with real coefficients.
In this isomorphism the action of ${\Gamma}$ translates in
$M_{2\times 2}(\mathbb R)$ as the linear action by conjugation of
$\phi_0({\Gamma})\subset O(2)$. Since $\phi_0({\Gamma})$ is
dihedral, $H^1(P,\mathbb R^2\times 0)^{{\Gamma}}\cong \mathbb R$.
\endproof

With a similar argument one can prove:

\begin{lem}\label{lem:H0xR} $H^*(\OO;0\times\mathbb R)\cong 0$.\qed
\end{lem}

\begin{cor} $H^*(M;0\times \mathbb R)\cong 0$.
\end{cor}

\proof We apply the Mayer-Vietoris exact sequence to the pair
$(\NN(\Sigma), M)$, where $\NN(\Sigma)$ is a tubular neighborhood of
$\Sigma$, so that $\NN(\Sigma)\cup M=\OO$ and $\NN(\Sigma)\cap
M\simeq T^2$. By Lemma~\ref{lem:H0xR}, we have an isomorphism:
\[
H^*(M;0\times \mathbb R)\oplus H^*(\NN(\Sigma);0\times \mathbb
R)\cong H^*(T^2;0\times \mathbb R).
\]
Since the meridian $m$ belongs to $\pi_1M$ and $\rho_0(m)$ acts
non-trivially on $0\times \mathbb R$, it follows that
$
    H^0(T^2;0\times \mathbb R)\cong H^0(\pi_1T^2;0\times \mathbb R)
    \cong (0\times \mathbb R)^{ \pi_1T^2}\cong 0.
$ By duality $H^2(T^2;0\times\mathbb R)\cong 0$, and by Euler
characteristic, $H^1(T^2;0\times\mathbb R)\cong 0$. \qed

\begin{lem} \label{lem:final} $H^1(M;su(2))\cong \mathbb
R$. In particular $H^1(M;\mathbb R^2\times 0)\cong \mathbb R$.
\end{lem}

\proof We apply a Mayer-Vietoris argument to the pair
$(M,\NN(\Sigma))$. Since $M\cup\NN(\Sigma)=\OO$ and
$M\cap\NN(\Sigma)\simeq T^2$,
 we have an exact sequence:
\[
 H_1(T^2,su(2))
\mathop{\longrightarrow}^{i_1\oplus i_2} H_1(M,su(2)) \oplus
H_1(\NN(\Sigma),su(2))\mathop{\longrightarrow}^{j_1-j_2}
H_1(\OO,su(2))
\]
where $i_1$, $i_2$, $j_1$ and $j_2$ are the natural maps induced by
inclusion. Notice that $j_1\circ i_1=j_2\circ i_2$ by exactness. We
have divided the proof in several steps.
\begin{enumerate}[(1)]
\item \emph{$H_1(T^2,su(2))\cong \mathbb R^2$ and $\{d\alpha_l, d\alpha_m\}$ is a basis for
$H_1(T^2,su(2))$.}

 This follows from the local
properties of the variety of representations $R(T^2,SU(2))$. See
\cite{PoOne}, for instance.

\item \emph{$j_2i_2(d\alpha_l)=j_1i_1(d\alpha_l)\neq 0 $. In particular
it is a basis for $H_1(\mathcal O,su(2))$.}

The proof that $j_1i_1(d\alpha_l)\neq 0$ uses the same argument as
the proof of Lemma~\ref{lem:w}. More precisely, since
$\pi\circ\hol(l)$ is a nontrivial translation, the Kronecker pairing
between the cocycle $z_q=\TRANS_q\circ\pi\circ\hol$ and $d\alpha_l$
does not vanish (Prop.~9.6 from \cite{PoTwo}). Thus $d\alpha_l\neq 0
$ when viewed in  $H^1(\pi_1\mathcal O,su(2))$. Since $H^1(\mathcal
O,su(2))\cong H^1(\pi_1\mathcal O,su(2))\cong\mathbb R$, by
Lemmas~\ref{lem:HRx0} and \ref{lem:H0xR}, it is clear that this
element is a basis.

\item $i_2(d\alpha_m)=0$.

This follows easily from the computation of
$H_1(\NN(\Sigma),SU(2))$, because  $m$ has order two, and
therefore it is rigid (see \cite{PoOne} for details).

\item $i_1\co H_1(T^2,su(2))\to H_1(M,su(2))$ has rank
one.

Since this map is Poincar\'e dual to $H_1(M,\partial M,su(2))\to
H_1(T^2,su(2))$, this follows from the long exact sequence of the
pair $(M,\partial M)$ and Step~1.

\item $i_1(d\alpha_m)=0$.

The proof is by contradiction. Assume that $i_1(d\alpha_m)\neq 0$.
Then by Step~4, $i_1(d\alpha_l)=\lambda i_1(d\alpha_m)$ for some
$\lambda\in\mathbb R$. In addition, since $i_2(d\alpha_m)=0$:
\[
j_1i_1(d\alpha_l)=\lambda j_1i_1(d\alpha_m)=\lambda
j_2i_2(d\alpha_m)=0
\]
which contradicts Step~2.

\item $H_1(M,su(2))\cong \mathbb R$.

By the previous steps $i_1\oplus i_2$ has rank one. The map
$j_1-j_2$ has also rank one, because $H_1(\mathcal O,su(2))\cong
\mathbb R$ and $j_1-j_2$ is surjective by Step~2. A standard
computation shows that $\dim_{\R}(H_1(\NN(\Sigma),SU(2)))=1$.
Therefore $H_1(M,su(2))\cong \mathbb R$.
\end{enumerate}
This finishes the proof of the lemma.\qed

\medskip

{\bf Acknowledgement}\qua This research was partially supported by MCYT
through grant BFM2000--0007.


\begin{thebibliography}{[B]}

\let\olditem\bibitem
\def\bibitem#1]#2{\olditem{#2}}



\bibitem[BA]{Leila} \textbf{L Ben Abdelghani}, \emph{Espace des
    repr\'esentations du groupe d'un noeud dans un groupe de Lie}, Thesis
         U. de Bourgogne (1998)

\bibitem[Bre]{Bre} \textbf{G\,E Bredon}, \emph{Introduction to compact transformation groups},
 Pure and Applied Mathematics, Vol. 46.
Academic Press, New York--London (1972)

\bibitem[BP]{BP} \textbf{M Boileau}, {\bf J Porti}, \emph{Geometrization of
    3--orbifolds of cyclic type}, Ast\'erisque  272 (2001)

\bibitem[BLP1]{BLP0} \textbf{M Boileau}, {\bf B Leeb}, {\bf J Porti},
     \emph{Uniformization of small 3--orbifolds}, C. R. Acad. Sci. Paris
    S\'er. I Math. 332 (2001)  57--62

\bibitem[BLP2]{BLP1} \textbf{M Boileau}, {\bf B Leeb}, {\bf J Porti},
     \emph{Geometrization of 3--dimensional orbifolds, Part I: Geometry of cone manifolds},
      preprint (2002)

\bibitem[BLP3]{BLP2} \textbf{M Boileau}, {\bf B Leeb}, {\bf J Porti},
\emph{Geometrization of 3--dimensional orbifolds}, preprint (2002)

\bibitem[CEG]{CEG} \textbf{R\,D Canary}, {\bf D\,B\,A Epstein}, {\bf P
    Green}, \emph{Notes on notes of Thurston}, from: ``Analytical and
    Geometric Aspects of Hyperbolic Space" (D\,B\,A Epstein, editor)
    London Math. Soc. Lecture Notes Ser. 111 Cambridge
    Univ. Press, Cambridge (1987) 3--92

\bibitem[CHK]{CHK}
    \textbf{D Cooper}, {\bf C Hodgson}, {\bf S Kerchkoff},
    \emph{Three dimensional Orbifolds and Cone Manifolds},
      Mathematical Society of Japan Memoirs 5 (2000)

\bibitem[Cul]{Cul}
   \textbf{M Culler},
   \emph{Lifting representations to covering groups},
    Adv. in Math.
    {59} (1986) 64--70

\bibitem[CS]{CS}
  \textbf{M Culler}, {\bf P Shalen},
  \newblock \emph{Varieties or group representations and splittings of
  3--manifolds},
  \newblock  Ann. of Math.
  {117} (1984) 401--476

\bibitem[HPS]{HPS}
    \textbf{M Heusener}, {\bf J Porti}, {\bf E Su\'arez},
    \emph{Regenerating singular hyperbolic  structures from $\textrm{Sol}$},
    J. Differential Geom. 59 (2001) 439--478 

\bibitem[Hod]{Hod}
    \textbf{C Hodgson},
     \emph{Degeneration and
    Regeneration
     of Hyperbolic Structures on
     Three-Manifolds}, Thesis, Princeton University  (1986)

\bibitem[HK]{HK}
     \textbf{C Hodgson}, {\bf S Kerckhoff},
     \emph{Rigidity of
    hyperbolic
    cone-manifolds and hyperbolic Dehn Surgery}, J. Diff. Geom. 48
    (1998)
    1--59

\bibitem[Koj]{Kojima}
    \textbf{S Kojima},
     \emph{Deformations of hyperbolic
    3--cone-manifolds}, J. Diff. Geom.
    {49} (1998) 469--516

\bibitem[LM]{LM}
    \textbf{A Lubotzky}, {\bf A\,R Magid},
   \emph{Varieties of representations of finitely generated groups},
     Mem. of the Amer. Math. Soc.
     {58} (1985)


\bibitem[Po1]{PoOne}
    \textbf{J Porti},
    \emph{Torsion de Reidemeister pour les Vari\'et\'es
    Hyperboliques},
    Mem. Amer. Math. Soc.
    {128} (1997)

\bibitem[Po2]{PoTwo}
    \textbf{J Porti},
   \emph{Regenerating hyperbolic and spherical
     cone structures from Euclidean ones},
    Topology
    {37}  (1998) 365--392

\bibitem[Sua]{Sua}
    \textbf{E Su\'arez},
 \emph{Poliedros de
    Dirichlet
     de 3--variedades c\'onicas y sus deformaciones}, Thesis, U.
    Complutense de
     Madrid (1998)

\bibitem[Sco]{Scott}
    \textbf{P Scott},
     \emph{The geometries of
    $3$--manifolds},
     Bull. London Math. Soc.
     {15} (1983) 401--487

\bibitem[Th1]{ThurstonNotes}
    \textbf{W\,P Thurston},
     \emph{The Geometry
     and Topology of 3--manifolds},  Princeton Math. Dept. (1979)

\bibitem[Th2]{Thurston}
    \textbf{W\,P Thurston},
                 \emph{Three-dimensional
                 geometry and topology}, Vol.  1.  (S Levy, editor)  Princeton
                 University Press, Princeton, NJ (1997)

\bibitem[Whi]{Whitney}
    \textbf{H Whitney},
     \emph{On singularities of mappings of Euclidean spaces.
I. Mappings of the plane into the plane},
     Ann. of Math.
     {62} (1955) 374--410

\end{thebibliography}
\end{document}